\documentclass[twoside,10.5pt]{article}

\usepackage{amssymb,amsmath,amsthm,float,graphicx,color,amsmath,cases,graphics,verbatim,mathrsfs}
\usepackage{enumerate,fullpage,natbib}

\usepackage[margin=1in]{geometry}

\definecolor{lightgray}{rgb}{0.95,0.95,0.95}

\newtheorem{theorem}{Theorem}
\newtheorem{lemma}{Lemma}
\newtheorem{proposition}{Proposition}

\def\[{\lf[} \def\]{\ri]}  \def\br{\right} \def\lf{\left} \def\ri{\right}

\def\old#1{}
\def\pf{{\smallskip\noindent\it Proof.  }}
\def\qed{{\hfill $\blacksquare$}}
\def\to{\rightarrow}

\def\E{{\bf E}} \def\P{{\bf P}}
 
\def\argmin{{\rm argmin}}

\def\l{\lambda}

\def\ek{ \mathcal{E}^t}

\def\a{\gamma}

\def\cV{\mathcal{V}}
\usepackage{multirow}

\author{Mengdi Wang\\ \small Department of Operations Research and Financial Engineering, Princeton University, Princeton, NJ\\  \small email: mengdiw@princeton.edu}
\title{Randomized Linear Programming Solves the Discounted Markov Decision Problem In Nearly-Linear (Sometimes Sublinear) Run Time}
\date{September 1, 2017}

\def\cA{\mathcal{A}}

\def\cO{\mathcal{O}}

\def\cF{\mathcal{F}}
\def\cM{\mathcal{M}}
\def\cS{\mathcal{S}}

\def \gk{\mathcal{G}^t}
   \def\cE{\mathcal{E}}

\usepackage{algpseudocode,algorithm}

\def\S{|\cS|}
\def\A{|\cA|}
\def\suma{\sum_{a\in \cA}}
\def\sumi{\sum_{i\in\cS}} 
\def\sumj {\sum_{j\in\cS}} 
\def\cF{\mathcal{F}}
\def\cP{\mathcal{P}}

 \def\g{\Delta}
 \def\sumia{\sumi\suma}
\def\KL{D_{KL}}
\def\tO{\tilde\cO}
\def\cU{\mathcal{U}_{\theta,\bq}}
\def\cG{\mathcal{G}}

\begin{document}

\maketitle

\begin{abstract}
We propose a novel randomized linear programming algorithm for approximating the optimal policy of the discounted Markov decision problem. 
By leveraging the value-policy duality and binary-tree data structures, the algorithm adaptively samples state-action-state transitions and makes exponentiated primal-dual updates. We show that it finds an $\epsilon$-optimal policy using nearly-linear run time in the worst case. When the Markov decision process is ergodic and specified in some special data formats, the algorithm finds an $\epsilon$-optimal policy using run time linear in the total number of state-action pairs, which is sublinear in the input size. These results provide a new venue and complexity benchmarks for solving stochastic dynamic programs. 
\end{abstract}


\def\bp{\mathbf{p}}
\def\bpi{\mathbf{\mu}}
\def\bmu{\mathbf{\pi}}
\def\bv{\mathbf{v}}
\def\br{\mathbf{r}}
\def\be{\mathbf{e}}
\def\bq{\mathbf{q}}
\def\bx{\mathbf{x}}

\section{Introduction}

Markov decision process (MDP) is a fundamental model for sequential decision-making problems in dynamic and random environments. It models a stochastic control process in which a planner aims to make a sequence of decisions as the state of the process evolves. MDP serves as the basic mathematical framework for dynamic programming, stochastic control and reinforcement learning. It is widely applied in engineering systems, artificial intelligence, e-commerce and finance.  

We focus on the Discounted Markov Decision Problem (DMDP) in which one aims to make an infinite sequence of decisions and optimize some cumulative sum of discounted rewards. An instance of the DMDP can be described by a tuple $\mathcal{M}=(\cS,\cA,\cP, \br,\gamma)$, where $\cS$ is a finite state space of size $\S$, $\cA$ is a finite action space of size $\A$, $\gamma\in(0,1)$ is a discount factor, $\cP$ is the collection of state-to-state transition probabilities $\cP=\{ p_{ij}(a) \mid  i,j\in\cS,a\in\cA\}$, $\br$ is the collection of state transitional rewards $\br=\{r_{ij}(a) \mid  i,j\in\cS,a\in\cA \}$ where $r_{ij}(a) \in [0, 1]$. We also denote by $\br_a\in\Re^{|\cS|}$ the vector of expected state-transition reward under action $a$, where $\br_{a,i} = \sumj p_{ij}(a)r_{ij}(a)$.
Suppose that the decision process is in state $i$, if action $a$ is selected, the process moves to a next state $j$ with probability $p_{ij}(a)$ and generates a reward $r_{ij}(a) $.  
The input size of the DMDP tuple $\mathcal{M}$ is $\cO(\S^2\A)$. 

Our goal is to find the best sequence of actions to choose at all possible states in order to maximize the expected cumulative reward. More precisely, we want to find a (stationary) policy that specifies which action to choose at each state. A stationary and randomized policy can be represented by a collection of probability distributions $\bmu = \{\bmu_i\}_{i\in\cS} $, where $\bmu_{i}:\cA\mapsto [0,1] $ is a vector of probability distribution over actions at state $i$.  
We denote by $P^{\bmu}$ the transition probability matrix of the DMDP under a fixed policy $\pi$, where $P^{\bmu}_{ij} = \suma \bmu_{i}(a) p_{ij}(a)$ for all $i,j\in\cS$.
The objective of the DMDP is to find an optimal policy $\bmu^*$ such that the infinite-horizon sum of discounted rewards is maximized regardless of the initial state $i_0$:
$$\max_{\bmu} \E^{\bmu}\[ \sum^{\infty}_{t=1} \gamma^t r_{i_{t} i_{t+1}}(a_t) \mid i_0 \],$$
where $\{i_0, a_0, i_1,a_1, \ldots,i_t,a_t,\ldots\}$ are state-action transitions generated by the Markov decision process under the fixed policy $\bmu$, and the expectation $ \E^{\bmu}\[\cdot\]$ is taken over the entire trajectory. 
In total there are $\A^{\S}$ distinct deterministic policies.  


Despite its strong power of modeling, MDP is generally considered as a difficult problem due to the {\it curse of dimensionality}, especially for problems with large state and action spaces. 
There have been tremendous efforts in analyzing the complexity of MDP and its solution methods. Most existing studies focus on deterministic methods that find the exact optimal policy. Due to the curse of dimensionality, finding the exact optimal policy is often prohibitively difficult, especially in large-scale applications such as computer games and robotics. We ask the following question:
\begin{center}
{\it Is there a way to trade the precision of the exact optimal policy for a better time complexity?} 
\end{center}
Motivated by this question, we are interested in developing randomized algorithms that can approximate the optimal policy efficiently. In particular, we are interested in reducing the complexity's dependence on $\S$ and $\A$ - sizes of the state and action spaces. Throughout this paper, we measure the run-time complexity of an algorithm in terms of the total number of arithmetic operations, which include query to a specific entry of the input, addition, subtraction, multiplication, division and exponentiation. We use $\cO(1)$ to denote some absolute constant number, and we use $\tO(1)$ to hide $\hbox{polylog}\lf(\S,\A,\frac1\epsilon,\frac1{1-\gamma}\ri)$ factors.


\subsection{Our Approach and Technical Novelties}
In this paper, we develop a randomized linear programming method for solving the DMDP. 
It can be viewed as a special stochastic primal-dual method that takes advantages of three features: (1) adaptive action sampling according to the current randomized policy; (2) multiplicative policy updates using information projection onto a specifically constructed constraint set; and (3) using binary-tree data structures to simulate state transitions and make policy updates in nearly constant time. 
Let us we outline the development of our method and its analysis.


\begin{enumerate}
\item Our starting point is to formulate the nonlinear Bellman equation for the DMDP into a stochastic saddle point problem (see Section 3). The primal and dual variables correspond to the value and the policy, respectively. Our saddle point formulation involves specially chosen constraints and a weight vector, which are crafted to incorporate structural information and prior knowledge (if any) about the DMDP, such as the discount factor, magnitudes of reward, and range of ergodic distributions. In particular, the dual constraint can be viewed as an information set that contains all possible randomized policies and facilitates $\tO(1)$-time projection with respect to the some variant of the relative entropy.

\item  To aid the algorithm design, we develop two programming techniques: (1) We show that by processing the input transition probabilities into binary trees (using $\tO(\S^2\A)$ time), one can sample a single state transition of the Markov decision process using $\tO(1)$ arithmetic operations; (2) We show that a randomized policy can be represented using tree data structures, such that each coordinate update can be made in $\tO(1)$ time and each random action can be sampled  in $\tO(1)$ time. These two sampling techniques enable us to develop randomized algorithms that simulate state transitions and make appropriate dual update in $\tO(1)$ time. 

\item Our randomized algorithm has two main components. First, it uses adaptive importance sampling of state-action-state triplets. In other words, it simulates the Markov decision process using the current dual variable as the control policy to balance the exploration-exploitation tradeoff in estimating the optimal policy. Second, our method iteratively makes exponentiated and re-weighted updates in the dual variable and projects it onto an information set with respect to a specific divergence function. Here the update rule and the divergence are designed jointly in a way such that each iteration takes $\tO(1)$ arithmetic operations.

\item Analyzing the convergence is complicated by the use of adaptive action sampling and weighted exponentiated updates, which leads to substantial noises with unbounded second moments. 
As a result, the classical primal-dual analysis by \cite{nemirovski2009robust,juditsky2011solving} no longer works. 
To tackle this difficulty, we develop an independent convergence proof by analyzing the stochastic improvement of a particular relative entropy. We obtain a finite-time duality gap bound that characterizes how much the complementarity condition of the Bellman linear program is violated by the output dual variable.

\item Another critical piece of our analysis is to study the relation between the duality gap of the iterate and the efficiency loss of the output randomized policy. We show that when the Markov decision process is sufficiently ergodic, the duality gap provides a sharp estimate for the value loss of the output approximate policy.

\item Finally, we provide a meta algorithm that performs multiple independent trials of the randomized primal-dual iteration to get a good policy with high probability. To achieve this goal, we develop a subroutine for policy evaluation and show that it computes an $\epsilon$-accurate value in $\tO(\frac1{\epsilon^2(1-\gamma)^2})$ run time. We prove that the meta algorithm is able to select the best policy out of many candidates with probability arbitrarily close to 1.
\end{enumerate}

\subsection{Main Results} 

We analyze the run-time complexity of the proposed randomized algorithm for obtaining an $\epsilon$-optimal policy, i.e., a policy that achieves $\epsilon$-optimal cumulative reward (to be specified in more details later).  The run-time complexity of an algorithm is measured by the number of arithmetic operations.
Our main results are summarized as follows:
\begin{enumerate}

\item We show that the proposed algorithm finds an $\epsilon$-optimal policy $\hat\pi$
with probability at least $1-\delta$ in run time
$$\tO\lf( \frac{\S^3 \A }{(1-\gamma)^6\epsilon^2}  {\log\lf(\frac1{\delta}\ri) } \ri).$$
We recall that the input size of the DMDP is $\cO(\S^2\A)$. This result establishes a {\it nearly-linear run time}. This is the first computational complexity result for using randomized linear programming to solve DMDP. 

\item Here comes the more interesting result: In the case where the decision process is ergodic under any stationary policy, we show that the algorithm finds an $\epsilon$-optimal policy with probability at least $1-\delta$ in run time
$$\tO\lf( \S^2\A + \frac{\S \A  }{(1-\gamma)^4\epsilon^2}  {\log\lf(\frac1{\delta}\ri)}\ri).$$
The first term $\tO\lf( \S^2\A\ri)$ is due to an initialization step that preprocesses the input data (consisting mainly of arrays of transition probabilities) into a tree-based sampler. This run time is {\it linear with respect to the input size}. Although requiring an additional ergodicity assumption, it has better dependence on $\S,\A$ than the best known simplex method  (also the policy iteration method) \cite{ye2011simplex},\cite{scherrer2013improved} and the value iteration method \cite{littman1995complexity}.

\item In addition, the preprocessing step can be omitted or expedited if the input data are given in a suitable data structure that can be directly used as a sampler (e.g., binary trees or arrays of cumulative sums). 
In these cases, the complexity upper bound reduces to 
$$\tO\lf( \frac{\S \A  }{(1-\gamma)^4\epsilon^2}{\log\lf(\frac1{\delta}\ri)} \ri) \ll \S^2\A,$$
where the inequality holds when $\S$ is sufficiently large.
This is a surprising {\it sublinear run-time complexity}.  In other words, it is possible to compute a near-optimal policy without even reading all entries of the input data. To the author's best knowledge, this is the first sublinear run-time result for DMDP.
\end{enumerate}

We compare the above complexity upper bounds with recent results on the computational complexity lower bound of DMDP \cite{chen2017lower}. It shows that any algorithm needs at least $\Omega(\S^2 \A)$ run time to get an $\epsilon$-approximate policy with high probability in general. It also shows that the lower bound reduces to  $\Omega(\frac{\S\A}{\epsilon})$ when the input data are in the format of binary trees or cumulative sums (for which the preprocessing step can be skipped), making sublinear-time algorithms possible. Comparing our main results with the lower bounds, we observe a counter-intuitive phenomenon:
{\it The computational complexity of DMDP depends on the input data structure.}



\paragraph{Notations} All vectors are considered as column vectors. For a vector $\bx\in\Re^n$, we denote by $x_i$ or $x(i)$ its $i$-th component, denote by $\bx^{\top}$ its transpose, and denote by $\|\bx\| =\sqrt{\bx^{\top}\bx}$ its Euclidean norm. 
We denote by $\be = (1,\ldots,1)^{\top}$ the vector with all entries equaling 1, and we denote by $\be_i$ the vector with its $i$-th entry equaling $1$ and other entries equaling $0$.
For a positive number $x$, we denote by $\log x$ the natural logarithm of $x$.
For two probability distributions $p,q$ over a finite set $X$, we denote by $\KL(p||q) $ their Kullback-Leibler divergence, i.e., $\KL(p||q) = \sum_{x\in X} p(x) \log \frac{p(x)}{q(x)}$.

\section{Related Literatures}

Three major approaches for solving MDP are the value iteration method, the policy iteration method, and linear programming methods. See the textbooks \citep{bertsekas1995dynamic, bertsekas1995neuro, puterman2014markov, bertsekas2013abstract} and references therein for more detailed surveys on MDP and its solution methods. 

Bellman \cite{bellman1957dynamic} developed the value iteration as a successive approximation method to solve the nonlinear fixed-point Bellman equation. Its convergence and complexity have been thoroughly analyzed; see e.g. \cite{tseng1990solving, littman1995complexity}.  The best known complexity result for value iteration is $\cO(\S^2\A L \frac{\log(1/(1-\gamma))}{1-\gamma})$, where $L$ is the number of bits to represent the input ($L\geq \S^2\A$). Value iteration can be also used to find an approximate solution in $\cO(\S^2\A \frac{\log(1/\epsilon(1-\gamma))}{1-\gamma})$.
Later it was shown by \cite{feinberg2014value} that the value iteration method is not strongly polynomial for DMDP.
Policy iteration was developed by Howard \cite{howard1960dynamic}, and its complexity has also been analyzed extensively; see e.g. \cite{mansour1999complexity,ye2011simplex,scherrer2013improved}. 
Ye \cite{ye2011simplex} showed that policy iteration (which is a variant of the general simplex method for linear programming) is strongly polynomial and terminates in $\cO(\frac{S^2A}{1-\gamma}\log(\frac{\S}{1-\gamma}))$ number of iterations. 
Later \cite{scherrer2013improved} improved the result by removing the $\log \S$ factor.
Not long after the development of value and policy iterations, \cite{d1963probabilistic} and \cite{de1960problemes} discovered that the Bellman equation can be formulated into an equivalent linear program. 
It followed that one can apply linear programming method such as the simplex method by Dantzig \cite{dantzig2016linear} to solve MDP {\it exactly}. Later \cite{ye2005new} designed a combinatorial interior-point algorithm (CIPA) that solves the DMDP in strongly polynomial time. 
Recent developments \cite{lee2014path, lee2015efficient}
showed that linear programs can be solved in $\tO(\sqrt{\hbox{rank}(A)} )$ number of linear system solves, which, applied to DMDP, leads to a run time of $\tO( \S^{2.5}\A L )$.
We also note that there have been many methods for approximate linear programming. However, they do not apply to DMDP directly because an $\epsilon$ error in the linear programming formation of the DMDP might lead to arbitrarily large policy error.


\begin{table}[h!]\footnotesize
\begin{center}
\begin{tabular}{|c|c|c|}
\hline
Value Iteration & $\S^2\A L \frac{\log(1/(1-\gamma))}{1-\gamma}$ and $\S^2\A \frac{\log(1/(1-\gamma)\epsilon)}{1-\gamma}$ & \cite{tseng1990solving, littman1995complexity}
\\ \hline
Policy Iteration (Block Simplex) &  $\frac{\S^4\A^2}{1-\gamma} \log(\frac1{1-\gamma})$ & \cite{ye2011simplex},\cite{scherrer2013improved}
\\  \hline
LP Algorithm  &   $\tO( \S^{2.5}\A L )$ &\cite{lee2014path}
\\  \hline
   Combinatorial Interior Point Algorithm   &    $\S^4\A^4\log\frac{\S}{1-\gamma}$ & \cite{ye2005new}
   \\  \hline
       {Phased Q-Learning}    &  $\tO(\S^2\A+\S^2\A+\frac{\S\A}{(1-\gamma)^6\epsilon^2})$  & \cite{kearns1999finite} and This Work
   \\  \hline
    {Randomized Primal-Dual Method}    & $\tO ( \frac{\S^3 \A}{(1-\gamma)^6\epsilon^2}) $   & Main Result 1
     \\  \hline
    {Randomized Primal-Dual Method}    & $\tO (\S^2\A + \frac{\S \A}{(1-\gamma)^4\epsilon^2}) $ under ergodicity assumption  & Main Result 2
     \\  \hline
        {Randomized Primal-Dual Method}    & $\tO (\frac{\S \A}{(1-\gamma)^4\epsilon^2}) $ under ergodicity assumption and special input format & Main Result 3\\
          \hline
        {Lower Bound}    & $\Omega(\S^2\A) $  & \cite{chen2017lower}\\
                 \hline
        {Lower Bound}    & $\Omega\lf(\frac{\S\A}{\epsilon}\ri) $ under special input format & \cite{chen2017lower}\\
\hline
\end{tabular}
\end{center}
\caption{Run-Time Complexity for solving DMDP, where $\S$ is the number of states, $\A$ is the number of actions per state, $\gamma \in(0,1)$ is the discount factor, and $L$ is the total bit size to present the DMDP input.}
\label{default}
\end{table}

 A most popular method in reinforcement learning is the Q-learning method. It refers to a class of sampling-based variants of the value iteration. The work \cite{kearns1999finite} proposed the phased Q-learning mehod and proved that it takes $\tO(\frac{\S\A}{\epsilon^2})$ sample transitions to compute an $\epsilon$-optimal policy, where the dependence on $\gamma$ is left unspecified. If we carry out a more careful analysis (which is not available in \cite{kearns1999finite}), we can see that the sample complexity of phased Q-learning is actually $\tO(\frac{\S\A}{(1-\gamma)^6\epsilon^2})$. No run-time analysis is given explicitly for phased Q-learning in \cite{kearns1999finite}. Using the techniques developed in this work, we can apply the binary-tree sampling techniques described in Section 4.1 and apply Prop.\ \ref{prop-o1} to the phased Q-learning. Then we can implement the method using appropriate preprocessing and $\tO(1)$ run time per sample/update, leading to a total run time $\tO(\S^2\A+\frac{\S\A}{(1-\gamma)^6\epsilon^2})$. Such run-time analysis is not available in any existing literature. Unlike Q-learning, most reinforcement learning algorithms need far more than $\tO(1)$ arithmetic operations to process a single sample transition, so they are not efficient solvers for DMDP.



Table 1 summarizes the best-known run-time complexity, in terms of the total number of arithmetic operations, of solution methods for DMDP.  Our second result is the sharpest among existing methods, which however requires an additional assumption on the ergodicity of the Markov chains (which we believe to be mild). The upper bound result $\tO (\S^2\A + \frac{\S \A}{(1-\gamma)^4\epsilon^2}) $ nearly matches the lower bound $\Omega({\S^2\A})$ by \cite{chen2017lower}, except for the extra ergodicity assumption.
Our third main result reveals a surprising phenomenon: In the case where the input data are given in preprocessed data format that enables sampling, the complexity reduces to $\tO ( \frac{\S \A}{(1-\gamma)^4\epsilon^2})$, which is a sublinear run-time result and nearly match the lower bound $\Omega(\frac{\S\A}{\epsilon})$ for the same case by \cite{chen2017lower}.

Our new method and analysis is motivated by the line of developments on stochastic first-order optimization methods. They originate from the stochastic approximation method for the root finding problem; see \citep{KuY03, BMP90, Bor08}. They find wide applications in stochastic convex optimization, especially problems arising from machine learning applications. Several studies have been conducted to analyze the sample complexity of stochastic first-order methods, started by  the seminal work \cite{nemirovski2009robust} and followed by many others.
In particular, our proposed method is closely related to stochastic first-order method for saddle point problems. The earliest work of this type is by \cite{nemirovski2005efficient}, which studied a stochastic approximation method for convex-concave saddle point problem and gives explicit convergent rate estimates. Later the work \cite{juditsky2011solving} studied a class of stochastic mirror-prox methods and their rate of convergence for solving stochastic monotone variational inequalities, which contains the convex-concave saddle point as a special case.  
Another related work is by \cite{clarkson2012sublinear}, which developed a class of sublinear algorithms for minimax problems from machine learning. It showed that it is possible to find a pair of primal-dual solutions achieving $\epsilon$ duality gap in run time $\tO\lf(\frac1{\epsilon^2}\ri)$. 

Unfortunately, none of existing results on the general stochastic saddle point problem directly applies to the DMDP. There are several gaps to be filled. The first gap lies in the saddle point formulation of the Bellman equation. Prior to this work, it was not clear how to appropriately formulate the DMDP into a saddle point problem with appropriate constraints in order to maintain complexity guarantee and enable constant-time projection at the same time. The second gap lies in the implementation and run-time efficiency of algorithms.  Earlier works on stochastic mirror-prox methods mainly focused on the iteration/sample complexity.  In this work, we care about the overall run-time complexity for solving DMDP. To achieve the best run-time efficiency, we will provide an integrated algorithmic design that combines the mathematics together with programming techniques. The third gap lies in the proof of convergence. Our algorithm uses importance adaptive sampling of actions, which creates unbounded noises and disables the analysis used in \cite{nemirovski2005efficient, juditsky2011solving}. As a result, we have to come up with an independent primal-dual convergence analysis.
The fourth gap lies between the duality gap and the performance of the output policy. A small duality gap does not necessarily imply a small policy error. In this work, we aim to close all these gaps and develop efficient randomized algorithms with run-time guarantees.

Let us also mention that there are two prior attempts (by the author of this paper) to use primal-dual iteration for online policy estimation of MDP. 
The work \citep{CDC2016} and its journal version \cite{chen2016stochastic} considered a basic stochastic primal-dual iteration that uses Euclidean projection and uniform sampling of state-action pairs to solve DMDP and established a sample complexity upper bound $\cO(\S^{4.5}\A \epsilon^{-2}) $. No run-time complexity analysis is available.

\section{Bellman Equation, Linear Programs, and Stochastic Saddle Point Problem}

Consider a DMDP tuple $\mathcal{M}=(\cS,\cA,\cP, \br,\gamma)$.
For a fixed policy $\bmu$, the value vector $\bv^{\pi}\in\Re^{\S}$ is defined as
$$
v^{\pi}_{i} =\E^{\bmu}\[ \sum^{\infty}_{t=1} \gamma^t  r_{i_{t} i_{t+1} }(a_{t}) ~\Big\vert~ i_1= i\],\qquad \forall~i\in\cS,
$$
where $\E^{\bmu}\[ \cdot\]$ is taken over the state-action trajectory $\{i_1,a_1,i_2,a_2,\ldots\}$ generated by the Markov decision process under policy $\pi$.
The optimal value vector $\bv^*\in\Re^{\S}$ is defined as
\begin{align*}
v^*_{i} &= \max_{\bmu } \E^{\bmu}\[ \sum^{\infty}_{t=1} \gamma^t  r_{i_{t} i_{t+1}}(a_{t})~\Big\vert~i_1= i\]
=\E^{\bmu^*}\[ \sum^{\infty}_{t=1} \gamma^t r_{i_{t} i_{t+1}}(a_{t})~\Big\vert~i_1= i\],\qquad 
\forall~i\in\cS.
\end{align*}
According to the theory of dynamic programming  \cite{puterman2014markov, bertsekas1995dynamic}, a vector $\bv^*$ is the optimal value function to the DMDP if and only if it satisfies the following $|\cS|\times |\cS|$ system of equations, known as the {\it Bellman equation}, given by
\begin{equation*}\label{eq-bell}
\begin{split}
v^*_i  
= \max_{a\in\cA }&\Bigg\{  \gamma \sum_{j\in\cS} p_{ij}(a)v^*_j +  \sum_{j\in\cS}  p_{ij}(a)r_{ij}(a)\Bigg\},\qquad \forall~ i \in \cS,
\end{split} 
\end{equation*}
When $\gamma\in(0,1)$, the Bellman equation has a unique fixed point solution $\bv^*$, and it  equals to the optimal value vector of the DMDP. Moreover, a policy $\bmu^*$ is an optimal policy of the DMDP if and only if it attains the elementwise maximization in the Bellman equation. For finite-state DMDP, there always exists at least one optimal policy $\pi^*$. If the optimal policy is unique, it is also a deterministic policy. If there are multiple optimal policies, there exist infinitely many optimal randomized policies.

The nonlinear Bellman equation is equivalent to the following $|\cS|\times (|\cS| |\cA|)$ linear programming problem (see \cite{puterman2014markov} Section 6.9 and the paper \cite{de2003linear}):
\begin{equation}\label{eq-primal}\begin{split}
&\hbox{minimize  } (1-\gamma) \bq^{\top} \mathbf{v}\\
&\hbox{subject to} \left({I} - \gamma P_a\right)  \mathbf{v} - \mathbf{r}_a  \geq 0, \qquad \forall~a\in \cA,
\end{split}\end{equation}
where $\bq\in\Re^{\S}$ is an {\it arbitrary} vector of probability distribution satisfying $\be^{\top}\bq=1$ and $\bq>0$, $P_a\in\Re^{|\cS|\times |\cS|}$ is matrix whose $(i,j)$-th entry equals to $p_{ij}(a)$, ${I}$ is the identity matrix with dimension $|\cS|\times |\cS|$ and $\br_a\in\Re^{|\cS|}$ is the expected state-transition reward vector under action $a$, i.e., 
$r_a(i) = \sum_{j\in\cS} p_{ij}(a) r_{ij}(a)$ for all $ i\in\cS.$
The dual linear program of \eqref{eq-primal} is 
\begin{equation}\label{eq-dual}\begin{split}
&\hbox{maximize }  \sum_{a\in\cA} \mu_a^{\top} \mathbf{r_a}\\
&\hbox{subject to }  \sum_{a\in\cA} \left({I} - \gamma P_a^{\top}\right)  \mu_a = (1-\gamma) \mathbf{\bq}, ~~ \mu_a\geq 0,~~\forall~a\in\cA. 
\end{split}\end{equation}
It is well known that each deterministic policy corresponds to a basic feasible solution to the dual linear program \eqref{eq-dual}. A randomized policy is a mixture of deterministic policies, so it corresponds to a feasible solution of program \eqref{eq-dual}. We denote by $\mu^* = (\mu_a^*)_{a\in\cA}$ the optimal solution to the dual linear program \eqref{eq-dual}. If there is a unique dual solution, it must be a basic feasible solution. In this case, it is well known that the basis of $\mu^*$ corresponds to an optimal deterministic policy.



We formulate the linear programs \eqref{eq-primal}-\eqref{eq-dual} into an equivalent minimax problem, given by
\begin{equation}\label{eq-saddle}
\begin{split}
&\min_{\bv \in\mathcal{V} } \max_{\bpi\in \cU }  (1-\gamma) \bq^{\top} \bv + \suma \bpi_a^{\top} \lf( ( \gamma P_a - I ) \bv + \br_a\ri).
\end{split}
\end{equation}
We construct $\mathcal{V}$ and $ \cU$ to be the search spaces for the value and the policy, respectively, given by
\begin{align*}
\mathcal{V} 
&= 
\lf\{ \|\bv\|_{\infty}\leq \frac1{1-\gamma}, \bv \geq 0\ri\},\qquad 
 \cU = \lf\{ \be^{\top}\mu =1, \mu\geq 0, \suma \mu_a \geq \theta \bq \ri\},
\end{align*}
where $\theta$ is a small value to be specificied.  Let us verify that $\bv^*\in \cV$ and $\mu^*\in\cU$. Since all rewards $r_{ij}(a)$ belong to $[0,1]$, we can easily verify that $0\leq v^*_i\leq \frac1{1-\gamma}$ for all $i$, therefore $\bv^*\in \mathcal{V}$. By multiplying both sides of the dual constraint $\sum_{a\in\cA} \left({I} - \gamma P_a^{\top}\right)  \mu^*_a = (1-\gamma) \mathbf{\bq}$ with $\be^{\top}$, we also verify that $\be^{\top} \mu^* = 1$ because $\bq$ is a probability distribution.
In subsequent analysis, we will specify choices of $\theta,\bq$ so that $\suma \mu_a^* \geq \theta\bq$ and $\mu^*\in \cU$. As long as $\bv^*\in\cV$ and $\mu^*\in\cU$, we will be able to tackle the DMDP by approximately solving the saddle point problem \eqref{eq-saddle}.

There are many ways to formulate the minimax problem \eqref{eq-saddle} into a stochastic saddle point problem. A key observation is that the probability transition matrix $P_a$ and the reward vector $\br_a$ are naturally expectations of some random variables. For example, we can rewrite \eqref{eq-saddle} as follows
\begin{equation*}\label{eq-stoch-saddle}
\begin{split}
&\min_{\bv \in\mathcal{V} } \max_{\bpi\in \cU }  (1-\gamma) \bq^{\top} \bv + \suma \bpi_a^{\top} \lf( \lf( \gamma \sumi  \E_{j\mid i,a}\[ \be_i \be_j^{\top}\] - I \ri) \bv +  \sumi  \E_{j\mid i,a}\[ r_{ij}(a) \be_i\] \ri),
\end{split}
\end{equation*}
where the expectation $\E_{j\mid i,a}\[\cdot\]$ is taken over $j\sim p_{ij}(a)$ where $i$ and $a$ are fixed.
The preceding formulation suggests that one can simulate transitions $j\mid i,a$ of the Markov decision process in order to draw samples of the Langrangian function and its partial derivatives. In addition, we will show later that simulation of the Markov decision process is almost {\it free} in the sense that each transition can be simulated in $\tO(1)$ arithmetic operations. This motivates the use of a stochastic primal-dual iteration for solving the DMDP. 


\section{Algorithms}

We develop our main algorithms in this section. We first develop a few programming techniques that may be useful to all simulation-based methods for DMDP. Then we propose the randomized primal-dual algorithms and analyze their run-time complexity per iteration.

\def\bw{\mathbf{w}}

\subsection{Programming Techniques for Markov Decision Processes}

Randomized algorithms for DMDP inevitably involve simulating the Markov decision process and making policy updates. Our first step is to develop implmentation techniques for the two operations so that they take as little as $\tO(1)$ run time.


\begin{proposition}\label{prop-o1}
Suppose that we are given arrays of transition probabilities $\cP= (P_a)_{a\in \cA}$ and a randomized policy $\pi =\{\pi_{i,a}\}_{i\in\cS,a\in\cA}$. We are also given a stream of updates to the weight vectors, each update taking the form
$\pi_{i,a}\leftarrow \tilde \pi_{i,a}, \pi_i \leftarrow \pi_i/\|\pi_i\|_1$ for some $i\in\cS,a\in\cA$.
 There exists an algorithm that preprocesses $\cP$ in $\tO(\S^2\A)$ time, makes each update in $\cO(\log\A) $ time, and, for any initial state and at any time, generates a single state transition of the Markov decision process  under the current policy in $\cO(\log \A)+\cO(\log \S)$ time.
\end{proposition}
\pf
Let us apply the binary-tree scheme \cite{wong1980efficient}  to the computation of Markov decision process. We are given the transition probabilities $(P_a)_{a\in \cA}$ as arrays. For each state-action pair $(i,a)$, we preprocess each row vector of transition probabilities $P_a(i,\cdot)$ into a binary tree, where each leaf stores the value $P_a(i,j)$ for some $j$ and each node stores the sum of its two children. We need to construct $\S\A$ trees for all transition probabilities, and the preprocessing time is $\tO(\S^2\A)$. Then for a given state-action pair $(i,a)$, it is possible to generate a random coordinate $j$ with probability $p_{ij}(a)$ by drawing a random variable uniformly from $[0,1]$ and search for the leaf that corresponds to an interval that contains the random variable. This procedure takes $\cO(\log \S)$ arithmetic operations.

We represent a randomized policy $\pi$ using a collection of $\S$-dimensional vectors of nonnegative weights $\{\bw_{i}\}_{i\in\cS}$ such that $\pi_{i,a} = \frac{w_{i,a}}{\sum_{a\in\cA} w_{i,a}}$.  
Similar to $P_a(i,\cdot)$'s,  the weight vectors $\{\bw_{i}\}_{i\in\cS}$ can be represented as $\S$ binary trees, one for each $i\in\cS$. For any given $i\in\cS$, one can generate a random coordinate $a$ with probability $\pi_{i,a} = \frac{w_{i,a}}{\|\bw_i\|_1}$ using $\cO(\log \S)$ arithmetic operations. Now suppose that we want to make the policy update $\pi_{i,a}\leftarrow\tilde \pi_{i,a}, \pi_i \leftarrow \pi_i/\|\pi_i\|_1$. Then we update the corresponding leave from $w_{i,a}$ to $ \|\bw_{i}\|_1 \cdot \tilde \pi_{i,a} $, where $\|\bw_{i}\|_1$ is simply the value of the root. We also need to update the values of all nodes on the path from the root to the leaf $w_{i,a}$, so that each node remains the sum of its two children. This update takes $\cO(\log\A) $ operations.

Finally, suppose that we want to simulate the decision process and generate a state transition according to the current policy. For a given state $i\in\cS$, we first sample an action $a$ according to the policy and then sample a coordinate $j$ according to $P_a(i,\cdot),$ which takes $\cO(\log \A)+\cO(\log \S)=\tO(1)$ arithmetic operations in total.
\qed


\vspace{6pt}

Proposition \ref{prop-o1} implies that simulation-based methods for DMDP can be potentially very efficient, because the sampling and updating operations are computationally cheap. This result suggests an intriguing connection between the sample complexity for estimating the optimal policy and the run-time complexity for approximating the optimal policy.

\begin{algorithm}[h!]
\caption{Randomized Primal-Dual Method for DMDPs}\label{algo:primaldual}
\begin{algorithmic}[1]
\State {\bf Input:} DMDP tuple $\mathcal{M}=(\cS,\cA,\cP,\br,\gamma)$, $\bq = \frac1S \be,\theta = (1-\gamma)$, $T$.
\State Set $\bv = 0\in\Re^{\S}$ 
\State Set $\xi = \frac1{\S} \be \in\Re^{\S} $, $\pi_i =\frac1{\A}\be\in\Re^{\A}$ for all $i\in \cS$
\State  Set $\beta= (1-\gamma) \sqrt{\frac{\log \lf(\S\A +1\ri)}{2\S\A T}}, \alpha= \frac{\S}{2(1-\gamma)^2} \beta,  M =\frac1{1-\gamma}$
\State Preprocess the input probabilities $\cP=\{p_{ij}(a)\}$ for sampling
     \hfill   (Complexity $\tO(S^2A)$)
\For{$t = 1, 2, \ldots,T$}
  \State \label{step-sample-ia}Sample $i$ with probability $((1-\theta) \xi_i +\theta q_i)$
   \hfill   (Complexity $\tO(1)$) 
     \State \label{step-sample-ia}Sample $a$ with probability $\pi_{i,a}$
   \hfill   (Complexity $\tO(1)$) 
  \State \label{step-sample-j}Conditioned on $(i,a)$, sample $j$ with probability $p_{ij}(a)$
   \hfill   (Complexity $\tO(1)$)
  \State  Update the iterates by
     \hfill   (Complexity $\tO(1)$)
\begin{align*}
\g &\leftarrow   {\beta}  \cdot \frac{\gamma v_j - v_i +r_{ij}(a) - M }{ ((1-\theta) \xi_i +\theta q_i)\pi_{i,a}} \\
v_{i} & \leftarrow  \max\lf\{\min\lf\{v_{i} - {\alpha} 
\lf( \frac{(1-\gamma)q_{i} }{(1-\theta) \xi_i +\theta q_i} - 1\ri), \quad \frac{ 1 }{1-\gamma} \ri\},0\ri\}
\\
v_j &\leftarrow
\max\lf\{\min\lf\{v_j -  {\alpha}\gamma  , \quad \frac{1 }{1-\gamma} \ri\},0\ri\} 
\end{align*}
 \State  \label{step-dual}Update the iterates by
     \hfill   (Complexity $\tO(1)$)
\begin{align*}
\xi_i &\leftarrow \xi_i + \xi_i \pi_{i,a} \lf(\exp\lf\{ \g \ri\}-1\ri),
\xi \leftarrow \xi/\| \xi \|_1\\
\bmu_{i,a} &\leftarrow \bmu_{i,a} \cdot \exp\lf\{\g \ri\},
\bmu_i \leftarrow \bmu_i/\|\bmu_i\|_1
\end{align*}
\EndFor
\State\label{step-output} {\bf Ouput:} Averaged policy iterate $\hat \bmu_i=  {\frac1T\sum^{T}_{k=1} \bmu^t_{i} }$ for all $i\in\cS$
\end{algorithmic}
\end{algorithm}

\subsection{The Randomized Primal-Dual Algorithm}

Motivated by the saddle point formulation of Bellman's equation, we develop a randomized linear programming method to compute an approximately optimal policy. The method is given in Algorithm \ref{algo:primaldual}. It is essentially a randomized primal-dual iteration that makes updates to an explicit primal variable and an {\it implicit} dual variable. More specifically,
it makes iterative coordinate updates to three variables: $\pi$, $\bv$ and $\xi$. Here $\pi\in \Re^{\S\A}$ represents a randomized policy and is guaranteed to satisfy $\be^{\top} \pi_i =1$ and $\pi_i\geq 0$ for all $i\in\cS$ throughout the iterations. The vector $\bv \in\Re^S$ is the primal variable (also the value vector) and is guaranteed to satisfy $0\leq \bv\leq \frac1{1-\gamma} \be$ throughout. The vector $\xi\in{\Re^S}$ represents some distribution over the state space and satisfies $\be^{\top} \xi = 1,\xi\geq 0$ throughout. The policy $\pi$ and the distribution $\xi$ jointly give the dual variable $\mu$ according to
$$\bpi_{i,a} = \lf( (1-\theta) \xi_i + \theta  q_i \ri) \pi_{i,a},\qquad \forall~ i\in\cS,a\in\cA.$$
The dual variable $\mu$ does not appear in the iteration of Algorithm \ref{algo:primaldual}.  It is updated implicitly through updates on $\pi$ and $\xi$. We can verify that the implicit dual variable $\mu$ satisfies $\mu\in\cU$ throughout the iteration.

Also note that Algorithm \ref{algo:primaldual} uses adaptive importance sampling of state-to-state transitions according to the current policy. In particular, it samples a state-action pair $(i,a)$ with probability $((1-\theta) \xi_i +\theta q_i)\pi_{i,a}$ (Steps 7-8). The sampling distribution varies as the algorithm proceeds. Some state-action pairs will be sampled with higher and higher probability, meaning that the action becomes more favorable than other actions for the corresponding state. In the mean time, some state-action pairs will be sampled less and less frequently, while the corresponding probability $\pi_{i,a}$ reduces to 0 eventually. This can be viewed as a form of ``reinforcement learning", because the algorithm tend to sample more often the actions that empirically worked well. To account for the nonuniform sampling probability, each sample $\Delta$ is re-weighted by $\frac1{((1-\theta) \xi_i +\theta q_i)\pi_{i,a}}$.



Now we analyze the run-time complexity for each step of Algorithm \ref{algo:primaldual}.
Step 5 preprocesses the input probabilities $\cP=\{p_{ij}(a)\}$  into $\S\A$ tree samplers, one for each state-action pair. Each of the sampler corresponds to a probability vector of $\S$ dimension. The preprocessing time for each $m$-dimensional vector is $\tO(m).$ Therefore the total processing time is $\tO(\S^2\A)$. 
According to Proposition \ref{prop-o1}, Steps 7-9 
take $\tO(1)$ arithmetic operations each, provided that  one updates and sample from $\xi,\pi$ using the binary-tree data structures.
Step 10 updates two entries of the value vector $\bv$ and uses run time $\tO(1).$
According to Proposition \ref{prop-o1}, Step  \ref{step-dual} can also be implemented using $\tO(1)$ arithmetic operations.
Step 13 requires taking average of past iterates $\{\pi^t\}^T_{t=1}$. This can be achieved by maintaining an additional variable to record the running multiplicative weights and running average (e.g., using a special binary tree structure). The additional storage overhead is $\cO(\S\A)$, and the additional computation overhead is $\tO(1)$ per iteration.

To sum up, the total run-time complexity of Algorithm \ref{algo:primaldual} consists of two parts: complexity of preprocessing and complexity per iteration. Preprocessing takes $\tO(S^2A)$ run time, which can be skipped if the input data are given in some data structure that enables immediate sampling (e.g., sorted arrays \cite{bringmann2012efficient} or binary trees \cite{wong1980efficient}). 
Each iteration of Algorithm \ref{algo:primaldual} takes $\tO(1)$ arithmetic operations.  

\subsection{The Meta Algorithm}

Finally, we are ready to develop a meta algorithm that computes an approximately optimal policy with probability arbitrarily close to 1. The idea is to run Algorithm \ref{algo:primaldual} for a number of independent trials, perform approximate value evaluation to the output policies, and select the best out of the candidates.

\begin{algorithm}[h!]
\caption{Meta Algorithm}\label{algo:meta}
\begin{algorithmic}[1]
\State {\bf Input:} DMDP tuple $\mathcal{M}=(\cS,\cA,\cP,\br,\gamma)$, $\bq = \frac1S \be,\theta = (1-\gamma)$, $T$.
\State 
Run Algorithm \ref{algo:primaldual} for $K = \Omega(\log\frac1\delta)$ independent trials with precision parameter $\frac{\epsilon}2$ and obtain output policies 
$\bmu^{(1)},\ldots, \bmu^{(K)}$. 

\State  For each output policy $\pi^{(k)}$ and initial distribution $\bq$, conduct approximate value evaluation and obtain an approximate evaluation $\bar Y^{(k)}$ with precision level $\frac{\epsilon}{2}$ and fail probability $\frac{\delta}{2K}$. 

\State  Output $\hat\pi = \bmu^{(k^*)}$ such that $k^* = \hbox{argmax}_{k=1,\ldots,K} \bar Y^{(k)}$. 

\State\label{step-output} {\bf Ouput:} Averaged dual iterate $\hat \bmu_i=  {\frac1T\sum^{T}_{k=1} \bmu^t_{i} }$ for all $i\in\cS$
\end{algorithmic}
\end{algorithm}

In Algorithm \ref{algo:meta}, Step 3 approximately computes the cumulative discounted rewards for all candidate policies and a fixed initial state distribution $\bq$. The aim is to find an $\epsilon$-approximation to each cumulative reward with probability at least $1- \frac{\delta}{2K}$. Its implementation and complexity will be discussed and analyzed in Section 5.3.

\section{Complexity Analysis}

We develop our main results in a number of steps. We first study the convergence of some duality gap of the randomized primal-dual iteration given by Algorithm 1. Then we study how to quantify the quality of the output randomized policy from the duality gap upper bound. Third we show how to approximately evaluate a randomized policy in a small run time. Finally we connect all the dots and establish the overall run-time complexity of Algorithms 1 and 2.

\subsection{Duality Gap Analysis}


In this section, we study the convergence of Algorithm \ref{algo:primaldual}. We observe that each iteration of Algorithm \ref{algo:primaldual} essentially performs a primal-dual update. To see this, we define the auxiliary dual iterate $\bpi^t_{i,a}$ as
$$\bpi^t_{i,a} = \lf( (1-\theta) \xi^t_i + \theta  q_i \ri) \pi^t_{i,a},\qquad \forall~ i\in\cS,a\in\cA.$$
Informally speaking, we see that the $t$-th iteration of Algorithm \ref{algo:primaldual} takes the form 
\begin{equation}\label{eq-pditer}
\begin{split}
\bv^{t+1} &=\argmin_{ \bv\in \mathcal{V}}\lf\{  (\partial_{\bv} L(\bv,\mu^t) + \varepsilon^{t+1})^{\top} (\bv-\bv^t) + \frac1{2\beta} \|  \bv - \bv^t \|^2 \ri\},
\\
\mu^{t+1} &= \argmin_{\mu \in { \cU}} \lf\{   (\partial_{\mu} L(\bv^t,\mu) + \omega^{t+1} )^{\top} (\mu - \mu^t) + \frac1{2\alpha} \Phi(\mu; \mu^t)\ri\},
\end{split}
\end{equation}
where $\varepsilon^{t+1} $ and $\omega^{t+1} $ are zero-mean random noises due to the sampling step, $\Phi$ is a Bregman divergence function given by
\begin{align*}
\Phi(\mu; \hat \mu) 
&=  (1-\theta) \KL(\lambda || \hat \lambda ) + \theta \sumi   q_i \KL(\bmu_i || \hat \bmu_i ) ,\end{align*}
where $(\lambda,\bmu)$ is determined by $\mu$ and $(\hat\lambda,\hat\bmu)$ is determined by $\hat\mu$. 
The dual feasible region $\cU$ plays the role of an ``information set," in which we search for the optimal policy. The information set $\cU$ shall be constructed to characterize properties and additional prior knowledge (if there is any) regarding the DMDP.
Clearly, the set $\cU$ and the divergence function $\Phi$ are determined by the input parameters $\theta,\bq$. 
We will specify values of the parameters $\theta,\bq$ in subsequent analysis.

One might attempt to analyze iteration \eqref{eq-pditer} using the general analysis for stochastic mirror-prox iterations by \cite{nemirovski2005efficient, juditsky2011solving}.  Unfortunately, this would not work. The general primal-dual convergence analysis given by \cite{nemirovski2005efficient, juditsky2011solving} requires
$\E\[ \| \partial_{\mu} L(\bv^t,\mu) + \omega^{t+1} \|_*^2 \mid \cF_t\]<\sigma^2 $ for some appropriate norm $\|\cdot\|_*$  and a finite constant $\sigma$, where $\cF_t$ denotes the collection of random variables up to the $t$-th iteration. Our Algorithm \ref{algo:primaldual} cannot be treated in this way, because it samples the partial gradients using adaptive weights (Steps 7,8) and re-weighted samples (Step 10). In particular, for a given state $i$, the action $a$ is sampled with probability $\pi^t_{i,a}$ and the corresponding partial gradient $\Delta =   {\beta}  \cdot \frac{\gamma v_j - v_i +r_{ij}(a) - M }{ ((1-\theta) \xi_i +\theta q_i)\pi_{i,a}} $ has been re-weighted with $1/\pi^t_{i,a}$ to maintain unbiasedness, leading to large variances on the order of $\pi^t_{i,a}(1/\pi^t_{i,a})^2 = 1/\pi^t_{i,a}$. 
Therefore, each iteration Algorithm \ref{algo:primaldual} suffers from unbounded noises in the following sense
$$\E\[  \lf\| \partial_{\mu} L(\bv^t,\mu^t) + \omega^t \ri\|_{*}^2 \mid \cF_t\] \geq\Theta\lf( \min_{i\in\cS,a\in\cA}{\pi^t_{i,a}}(\frac1{\pi^t_{i,a}})^2\ri) \rightarrow \infty,\qquad \hbox{as}\ \  t\to\infty.$$
The lefthand side  can take arbitrarily large values and eventually go to infinity. This is because many $\pi^t_{i,a}$'s become increasingly close to zero, which is inevitable as the randomized policy $\pi^t$ converges to the optimal deterministic policy $\pi^*$. Due to this reason, the results and analyses given in \cite{nemirovski2005efficient, juditsky2011solving} do not apply.

To tackle this analytical difficulty, we develop an independent proof tailored to Algorithm \ref{algo:primaldual}. It involves analyzing the improvement of the relative entropy directly and constructing appropriate martingales. 
We obtain an finite-iteration upper bound on some expected ``duality gap".

\begin{proposition}[\bf Duality Gap Bound] \label{thm-dualitygap}
Let $\mathcal{M}=(\cS,\cA,\cP, \br,\gamma)$ be an arbitrary DMDP tuple, let $\bq\in\Re^{\S}$ be an arbitrary probability vector and let $\theta\geq 1-\gamma$.  Assume that the solution of the dual linear program \eqref{eq-dual} satisfies $\mu^*\in  \cU  $.
Let Algorithm \ref{algo:primaldual} iterate with the input $(\cM,\bq,\theta)$, then the sequence of iterates $\{(\xi^t,\pi^t)\}^T_{t=1}$ satisfies
$$\E\[ \sum_{a\in\cA,i\in\cS}\hat \bpi_{a,i} \lf( v_i^* - \a \sumj p_{ij}(a) v^*_j -\sumj p_{ij}(a) r_{ij}(a)  \ri)  \] \leq \frac{\sqrt{2\S (\A+1) \lf(\log(\S \A)+1\ri)}}{(1-\gamma) \sqrt{T}},$$
where $\hat \mu_{i,a} = \frac1T\sum^{T}_{t=1} ((1-\theta) \xi_i^t +\theta q_i)\pi_{i,a}^t$.
\end{proposition}

Proposition \ref{thm-dualitygap} establishes an expected upper bound on a nonnegative quantity that involves only the averaged dual variable $\hat\mu$ but not the primal variable $\bv$. It can be essentially viewed as a weighted sum of errors, which characterizes how much the following linear complementarity condition
$$\bpi_{a,i} \lf( v_i^* - \a \sumj p_{ij}(a) v^*_j -\sumj p_{ij}(a) r_{ij}(a)  \ri)  = 0,\qquad \forall a\in\cA,i\in\cS,$$
is violated for the linear programs \eqref{eq-primal}-\eqref{eq-dual}. Although we refer to this error quantity intuitively as a duality gap, it is actually not the typical minimax duality gap for saddle point problems (which was analyzed in \cite{nemirovski2005efficient, juditsky2011solving}).
We defer the detailed proof of Proposition \ref{thm-dualitygap} to Section 7.


\subsection{From Dual Variable To Approximate Policy}
%



In the following, we show that the duality gap of $\hat\mu$ gives an upper bound on the efficiency loss of the randomized policy $\hat \bmu$.

\begin{proposition}\label{prop-gaptoinf}
Let $\mathcal{M}=(\cS,\cA,\cP, \br,\gamma)$ be an arbitrary DMDP tuple, and let
$\bq=\frac1{\S}\be, \theta = (1-\gamma)$.
For any vector $\hat \mu\in \cU=  \lf\{ \be^{\top}\mu =1, \mu\geq 0, \suma \mu_a \geq \theta \bq \ri\}$, we let $\hat\pi$ be the corresponding randomized policy satisfying $\hat \pi_{i,a} = \frac{\hat\mu_{i,a}}{\suma \hat\mu_{i,a}}$ for all $i\in\cS,a\in\cA$. Then 
$$ \|\bv^* - \bv^{\hat\pi}\|_{\infty} \leq \frac{\S}{(1-\gamma)^2}
\sum_{a\in\cA,i\in\cS}\hat \bpi_{a,i} \lf( v_i^* - \a \sumj p_{ij}(a) v^*_j -\sumj p_{ij}(a) r_{ij}(a) \ri).
$$
\end{proposition}

\def\l{\lambda}
\pf
We denote for short that $Gap = \sum_{a\in\cA,i\in\cS}\hat \bpi_{a,i} \lf( v_i^* - \a \sumj p_{ij}(a) v^*_j -\sumj p_{ij}(a) r_{ij}(a) \ri)$. 
Using the fact $\hat\mu\in\cU$, we have $\sum_{a\in\cA}\hat \mu_{i,a} \geq \theta\bq = (1-\gamma)\bq. $ Then we have
\begin{align*}
Gap
&= \sum_{a\in\cA,i\in\cS}\hat \bpi_{i,a}  ( \bv^* - \a P_{a} \bv^* - \br_{a} )_i \\
&\geq (1-\gamma) \sum_{i\in\cS} q_i  \suma \hat\bmu_{a,i}( \bv^* - \a P_{a} \bv^* - \br_{a} )_i \\
&= (1-\gamma)  \sum_{i\in\cS}  q_i  ( \bv^* - \a P^{\hat\pi} \bv^* - \br^{\hat\pi} )_i ,
\end{align*}
where $\br^{\hat\pi}$ denotes the vector with $r^{\hat\pi}_i= \suma \hat\pi_i(a) \sumj p_{ij}(a) r_{ij}(a)$. We note that the Bellman equation for a fixed policy $\hat\pi$ is given by 
$\bv^{\hat\pi} =\gamma P^{\hat\pi} \bv^{\hat\pi}+\br^{\hat\pi}$.
Because $\bv^*$ is the optimal value vector of the DMDP, we have $( \bv^* - \a P_{a} \bv^* - \br_{a} )_i \geq 0$ for all $i\in\cS$. It follows that
$( \bv^* - \a P^{\hat\pi} \bv^* - \br^{\hat\pi} )_i  \geq 0$ for $i\in\cS$, therefore
$$0\leq \bv^* - \a P^{\hat\pi} \bv^*- \br^{\hat\pi} 
= \bv^* - \a P^{\hat\pi} \bv^*- (\bv^{\hat \bmu }  - \a P^{\hat\pi} \bv^{\hat\pi})
= (I-\gamma P^{\hat\pi}) (\bv^*- \bv^{\hat\pi}) .
$$
Using the fact $\bq=\frac1{\S} \be$, we further have
$
Gap 
=  (1-\gamma) (\bq)^{\top} (I-\gamma P^{\hat\pi}) (\bv^*- \bv^{\hat\pi})
\geq \frac{1-\gamma}{\S} \|  (I-\gamma P^{\hat\pi}) (\bv^*- \bv^{\hat\pi})\|_{\infty}.
$
We use the triangle inequality and the matrix norm inequality $\|A x \|_{\infty} \leq \|A \|_{\infty} \|x \|_{\infty} $ to obtain
$\|  (I-\gamma P^{\hat\pi}) (\bv^*- \bv^{\hat\pi})\|_{\infty} \geq \|   \bv^*- \bv^{\hat\pi}\|_{\infty}
 -  \|  \gamma P^{\hat\pi}(\bv^*- \bv^{\hat\pi})\|_{\infty}
 \geq \|   \bv^*- \bv^{\hat\pi}\|_{\infty}
 -  \|  \gamma P^{\hat\pi}\|_{\infty} \|(\bv^*- \bv^{\hat\pi})\|_{\infty}
  = (1-\gamma )\| \bv^*-\bv^{\hat\pi}\|_{\infty}.$
 It follows that $Gap  \geq \frac{(1-\gamma )^2}{\S} \| \bv^*-\bv^{\hat\pi}\|_{\infty}.$
 \qed
 
 \vspace{6pt}
 
 Next we will see that, the duality gap provides a sharper upper bound for the policy error when the associated Markov process is ``better" behaved.
For an arbitrary policy $\pi$, we define $\nu^{\pi}$ to be the {\it stationary distribution under policy $\pi$}, i.e.,
$\nu^{\pi}= (P^{\pi})^{\top}\nu^{\pi}$. 

 \begin{proposition}\label{prop-gaptoav}
Suppose that the Markov decision process specified by $\mathcal{M}=(\cS,\cA,\cP, \br,\gamma)$ is ergodic in the sense that $c_1 \bq\leq \nu^{\bmu} \leq c_2 \bq$ for some distribution vector $\bq$ and any policy $\pi$. Let $\hat \mu\in\cU$ where $\theta = 1-\gamma+\gamma \frac{c_1}{c_2}$, and let $\hat \pi_{i,a} = \frac{\hat\mu_{i,a}}{\suma \hat\mu_{i,a}}$. 
Then 
$$  \bq^{\top} \bv^* - \bq^{\top} \bv^{\hat\pi}   \leq \frac{c_2^2}{(1-\gamma)c_1^2}
\sum_{a\in\cA,i\in\cS}\hat \bpi_{a,i} \lf( v_i^* - \a \sumj p_{ij}(a) v^*_j -\sumj p_{ij}(a) r_{ij}(a) \ri).
$$
\end{proposition}
 \pf
 We have
\begin{align*}
\sum_{a\in\cA,i\in\cS}\hat \bpi_{a,i} ( \bv^* - \a P_{a} \bv^* - \br_{a} )_i 
&\geq \lf(1-\gamma + \gamma \frac{c_1}{c_2} \ri)\sum_{i\in\cS} q_i  \suma \hat\bmu_{i,a}( \bv^* - \a P_{a} \bv^* - \br_{a} )_i  \\
&= \lf(1-\gamma + \gamma \frac{c_1}{c_2} \ri) \sum_{i\in\cS}  q_i  
( \bv^* - \a P^{\hat\pi} \bv^* - \br^{\hat\pi} )_i \\
&\geq \lf(1-\gamma + \gamma \frac{c_1}{c_2} \ri) \sum_{i\in\cS}  \frac1{c_2} \nu^{\hat\pi}_i 
( \bv^* - \a P^{\hat\pi} \bv^* - \br^{\hat\pi} )_i  \\
&= \lf(1-\gamma + \gamma \frac{c_1}{c_2} \ri) \frac1{c_2}  \lf(\nu^{\hat\pi}\ri)^{\top}  (I - \a P^{\hat\pi}) (\bv^* - \bv^{\hat\pi}   ) \\
&= \lf(1-\gamma + \gamma \frac{c_1}{c_2} \ri)\frac1{c_2} (1-\gamma ) \lf( \nu^{\hat\pi}\ri)^{\top}  (\bv^*-\bv^{\hat\pi}) \\
&\geq  \frac{c_1^2}{c_2^2} (1-\gamma )  \bq^{\top} (\bv^*-\bv^{\hat\pi}) ,
\end{align*}
where the first inequality uses $\hat\mu \in\cU$, the second and third inequalities use $\nu^{\hat \pi}\in\cU$.
\qed

Prop.\ \ref{prop-gaptoav} suggests that, when all ergodic distributions of the Markov decision process under stationary policies belong to a certain range,  we get a sharper bound for the cumulative value of the output policy from the duality gap bound. 

\subsection{Approximate Policy Evaluation}

In Algorithm \ref{algo:meta}, we run Algorithm \ref{algo:primaldual} for multiple independent trials and identify the most successful one, in order to boost the probability of finding a good policy to be arbitrarily close to 1. 
Before we can do that, we need to be able to evaluate multiple candidate policies and select the best one out of many with high probability (Step 3 of Algorithm \ref{algo:meta}). 
In fact, we show that it is possible to approximately evaluate any policy $\pi$ within $\epsilon$ precision in run time $\tO(\frac1{\epsilon^2(1-\gamma)^2})$ for a pre-specified initial distribution.

\begin{proposition}[\bf Approximate Policy Evaluation]\label{lemma-pe}
Suppose we are given a DMDP tuple $\mathcal{M}=(\cS,\cA,\cP, \br,\gamma)$, a fixed randomized policy $\pi$, and an initial distribution $\bq$. Suppose that a state transition of the DMDP under $\pi$ can be sampled in $\tO(1)$ time,
there exists an algorithm that outputs an approximate value $\bar Y$ such that $\bq^{\top} \bv^{\pi}  -  {\epsilon} \leq \bar Y \leq  \bq^{\top} \bv^{\pi} $ with probability at least $1-\delta$ in run time 
$\tO(\frac1{\epsilon^2(1-\gamma)^2}\log(\frac1{\delta})) $.
\end{proposition}

\pf The approximate policy evaluation algorithm runs as follows:
For a given policy $\bmu$, we simulate the Markov decision process under policy $\bmu$ from the initial distribution $\bq$ for $n$ transitions and calculate the cumulative discounted reward $Y$. We repeat the simulation for $K$ independent times and return the average cumulative reward $\bar Y = \frac1K(Y_1+\cdots+Y_K)$. 

We observe that the unknown value $\bq^{\top} \bv^{\pi}$ is the expectation of the infinite cumulative discounted reward.
We pick $n$ sufficiently large such that expected $n$-period cumulative discounted reward is sufficiently close to  $\bq^{\top} \bv^{\pi}$. Since $r_{ij}(a)\in[\frac12,1]$ for $i,j,a,$, 
the cumulative discounted reward starting from the $(n+1)$th period is bounded by $\sum^{\infty}_{t=n} \gamma^t = \frac{\gamma^{n}}{1-\gamma}$ with probability 1. In particular, we pick $n$ such that $\frac{\epsilon}2 = \frac{\gamma^{n}}{1-\gamma}$, which suggests that $n= \lf(\log_\gamma (\frac{\epsilon(1-\gamma)}{2})\ri) =\tO(1)$. Therefore we obtain
$$  \bq^{\top} \bv^{\pi}  -  \frac{\epsilon}2 \leq \E\[Y_1\]  \leq  \bq^{\top} \bv^{\pi}  .$$
Note that $Y_1,\ldots,Y_K$ are i.i.d.\ random variables and $Y_k \in [0, \frac1{1-\gamma}]$ for all $k$.
By using the Azuma-Hoeffding inequality, we obtain that $\bar Y =  \frac1K\sum^{K}_{t=1} Y_t$ satisfies for any $\varepsilon>0$ that
\begin{equation*}
\P\lf(| \bar Y-\E\[ Y_1 \]| \geq  \varepsilon \ri) \leq 2 \exp\lf(-\frac{\varepsilon^2 K(1-\gamma)^2}2\ri).
\end{equation*}
By letting $\varepsilon=\frac{\epsilon}2$ and $K = \cO(\frac{1}{\epsilon^2(1-\gamma)^2} \log(1/\delta))$, we obtain that $|\bar Y- \E\[Y_1\]|\leq \frac{\epsilon}2$ with probability at least $1-\delta$.
It follows that 
$$\mathbf{P}\lf( \bq^{\top} \bv^{\pi}  -  {\epsilon} \leq \bar Y \leq  \bq^{\top} \bv^{\pi}  \ri) \geq 1-\delta
. $$
The number of sample state transitions is $K\cdot n = \tO(\frac1{\epsilon^2(1-\gamma)^2} \log\lf(\frac1\delta\ri)) $, which equals to the total run time.
\qed

\vspace{3pt}

\subsection{Run-Time Complexity For Algorithms \ref{algo:primaldual} and \ref{algo:meta}}



Finally, we are ready to develop the main results of this paper. 
Our first main result is given in Theorem 1. It establishes the run-time complexity for computing an $\epsilon$-optimal policy for arbitrary DMDP using the randomized linear programming methods given by Algorithms \ref{algo:primaldual} and \ref{algo:meta}.  

\begin{theorem}[\bf Run-Time Complexity for Arbitrary DMDP]\label{thm-1}

Let $\mathcal{M}=(\cS,\cA,\cP, \br,\gamma)$ be an arbitrary DMDP. For any $\epsilon \in(0,1)$, $\delta\in(0,1)$, let $\bq=\frac1{\S}\be, \theta = (1-\gamma)$, $ T = \Omega\lf(\frac{\S^3 \A \log(\S \A) }{(1-\gamma)^6\epsilon^2}   \ri).$ Then: 
\begin{enumerate}[(i)]
\item Algorithm \ref{algo:primaldual} outputs a policy $\hat\bmu$ satisfying $ \bv^{\hat \bmu}(i) \geq  \bv^*(i) - \epsilon$ for all $i\in\cS$ with probability at least $2/3$.
\item
Algorithm \ref{algo:meta} outputs a policy $\hat\bmu$ such that $ \bq^{\top }\bv^{\hat \bmu} \geq  \bq^{\top }\bv^* -\epsilon $ in run time
$$  \tO\lf(\S^2\A + \frac{\S^3 \A }{(1-\gamma)^6\epsilon^2} \log\lf(\frac1{\delta}\ri) + \frac{1}{\epsilon^2(1-\gamma)^2}\lf(\log\frac1\delta\ri)^2 \ri)$$
 with probability at least $1-\delta$.
 \end{enumerate}
\end{theorem}
\pf 
(i) We first show that the dual optimal solution $\mu^*$ to the linear program \eqref{eq-dual} satisfies $\mu^*\in \cU$ when $\theta= 1-\gamma,$ $\bq=\frac1S\be.$
To see this, we note the dual feasibility constraint 
$\suma (I-\gamma P_a)^{\top}\bpi^*_a = (1-\gamma)\bq,$
which implies a lower bound to the dual variable $\bpi$, given by
$\suma \bpi^*_a \geq (1-\gamma) \bq.$
Therefore $\mu^*\in \cU$ and the assumption of Proposition \ref{thm-dualitygap} is satisfied. 

We denote for short that $Gap = \sum_{a\in\cA,i\in\cS}\hat \bpi_{a,i} \lf( v_i^* - \a \sumj p_{ij}(a) v^*_j -\sumj p_{ij}(a) r_{ij}(a) \ri)$. 
Now we apply Prop.\ \ref{prop-gaptoinf} and the Markov inequality to obtain 
$\| \bv^*-\bv^{\hat\pi}\|_{\infty} 
\leq \frac{\S}{(1-\gamma )^2} Gap \leq \cO(\frac{S}{(1-\gamma )^2}) \E\[Gap\]$ with probability at least $2/3$.
 Therefore $\bv^\pi \geq \bv^* - \epsilon\be  $ 
 with probability $2/3$ when $\E\[Gap\] \leq \cO(\frac{\epsilon(1-\gamma)^2}{\S})$, which holds  if we let $T = \Omega\lf(\frac{\S^3 \A \log(\S \A) }{(1-\gamma)^6\epsilon^2}   \ri)$ and apply Prop.\ \ref{thm-dualitygap}.
 
 \vspace{4pt}
 
\noindent (ii)
Let us analyze Algorithm \ref{algo:meta} step by step. 
\begin{enumerate}
\item In Step 2 of Algorithm \ref{algo:meta}, it runs Algorithm \ref{algo:primaldual} for $K$ independent trials with precision parameter $\frac{\epsilon}2$ and generates output policies
$\bmu^{(1)},\ldots, \bmu^{(K)}$. The total run time is $\cO(\S^2\A)+ K\cdot T_{\frac{\epsilon}2}$, where $T_{\frac{\epsilon}2}$ is the time complexity for each run of Algorithm \ref{algo:primaldual}. According to (i), each trial generates an $\epsilon/2$-optimal policy with probability at least $2/3$.

\item  In Step 3 of Algorithm \ref{algo:meta}, for each output policy $\pi^{(k)}$, we conduct approximate value evaluation and obtain an approximate evaluation $\bar Y^{(k)}$ with precision level $\frac{\epsilon}{2}$ and fail probability $\frac{\delta}{2K}$. According to Lemma \ref{lemma-pe}, we have
$$  \bar Y^{(k)} - \bq^{\top} \bv^{\pi^{(k)}}\in\ [- \frac{\epsilon}2  ,0] ,
$$ 
with probability at least $1-\frac{\delta}{2 K}$. This step takes $K\cdot  \tO(\frac1{\epsilon^2(1-\gamma)^2} \log\lf(\frac{K}\delta\ri)) $ run time.

\item Step 4 of Algorithm \ref{algo:meta} outputs $\hat\pi = \bmu^{(k^*)}$ such that $k^* = \hbox{argmax}_{k=1,\ldots,K} \bar Y^{(k)}$. This step takes $\cO( K)$ run time.
\end{enumerate}

Now we verify that $\hat\pi$ is indeed near-optimal with probability at least $1-\delta$ when $K$ is chosen appropriately.
Let $\mathcal{K} =\lf\{k\mid \bq^{\top} \bv^{\pi^{(k)}} \geq \bq^{\top}\bv^{*} - \frac{\epsilon}2 \ri\} $, which can be interpreted as indices of successful trails of Algorithm \ref{algo:primaldual}. 
Consider the event where $\mathcal{K}\neq\emptyset$ and all policy evaluation errors belong to the small interval $[-\frac{\epsilon}2,0]$. 
In this case, we have $\bar Y^{(k)} \leq \bq^{\top} \bv^{\pi^{(k)}} $ for all $k$ and $\bar Y^{(k)} \geq   \bq^{\top} \bv^{\pi^{(k)}} - \frac{\epsilon}2 \geq \bq^{\top}\bv^{*}- {\epsilon}  $ if $k\in\mathcal{K}$. As a result, the output policy $\hat\pi$ which has the largest value of $\bar Y^{(k)} $ must be $\epsilon$-optimal.
We use the union bound to obtain 
\begin{align*}
\mathbf{P}\lf( \bq^{\top}\bv^{\hat\pi} <    \bq^{\top}\bv^{*} - \epsilon\ri)
&\leq  \mathbf{P}\lf(\lf\{ \mathcal{K} =\emptyset \ri\} \cup \lf\{ \exists k: \bar Y^{(k)} - \bq^{\top} \bv^{\pi^{(k)}}\notin [- \frac{\epsilon}2 ,0] \ri\}\ri)
\\
&\leq   \mathbf{P} \lf( \mathcal{K} =\emptyset \ri) +
\mathbf{P}\lf(  \exists k: \bar Y^{(k)} - \bq^{\top} \bv^{\pi^{(k)}}\notin [- \frac{\epsilon}2  ,0] \ri)
\\
&\leq \prod_{k=1}^K \mathbf{P} \lf(  \bq^{\top} \bv^{\pi^{(k)}} <  \bq^{\top}\bv^{*} -  \frac{\epsilon}2\ri)
+
\sum_{k=1}^K \mathbf{P}\lf(\bar Y^{(k)} - \bq^{\top} \bv^{\pi^{(k)}}\notin [- \frac{\epsilon}2  ,0] \ri)
\\
&\leq (1/3)^K +  K\cdot \frac{\delta}{2K} .
\end{align*}
By choosing $K  = \Omega (\log(1/\delta))$, we obtain $\mathbf{P}\lf( \bq^{\top}\bv^{\hat\pi} <    \bq^{\top}\bv^{*}  - \epsilon\ri)\leq \delta$. Then the output policy is $\epsilon$-optimal when initiated at distribution $\bq$ with probability at least $1-\delta$. 

According to Section 4, preprocessing of Algorithms \ref{algo:primaldual}-\ref{algo:meta} takes $\tO(\S^2\A)$ arithmetic operations and each iteration takes $\tO(1)$ arithmetic operations. The total run time is $\tO(\S^2\A + T_{\frac{\epsilon}2} \log\frac1\delta + \frac1{\epsilon^2(1-\gamma)^2} (\log\frac1\delta )^2)$.
\qed

\vspace{6pt}

Theorem \ref{thm-1} establishes a worst-case complexity for Algorithms \ref{algo:primaldual} and \ref{algo:meta}. The results apply to all instances of DMDP with $\S$ states, $\A$ actions per state, and a fixed discount factor $\gamma$. Theorem \ref{thm-1} does not require any additional property such as irreducibility or aperiodicity of the associated Markov chains. In fact, the results even hold for problems with transient states and/or multiple optimal policies. 

The cubic dependence $\S^3$ in Theorem \ref{thm-1} is not very satisfying. One factor of $\S$ comes from the duality gap bound (Prop.\ \ref{thm-dualitygap}), and two more factors come from rounding the duality gap to the policy error $\|\bv^{\hat\pi}  - \bv^*\|_{\infty}$ (Prop.\ \ref{prop-gaptoinf}). We conjecture that the term $\S^3$ can be improved to $\S^2$ (or even $\S$) using an improved algorithm and analysis. In addition, we conjecture that the complexity result should have a better dependence on $\frac1{1-\gamma}$, especially when all $P^{\pi}$'s have relatively large spectral gaps. These two questions are left open for future investigation.



\vspace{6pt}

Next we will see that, the DMDP becomes easier to solve when the associated Markov process is ergodic.
In what follows, we focus on the class of Markov decision processes where every stationary policy generates an ergodic Markov chain. 
For an arbitrary policy $\pi$, we define $\nu^{\pi}$ to be the {\it stationary distribution under policy $\pi$}, i.e.,
$\nu^{\pi}= (P^{\pi})^{\top}\nu^{\pi}$.
Our next main result shows that Algorithms \ref{algo:primaldual}-\ref{algo:meta} have significantly improved complexity for ergodic DMDP.

\begin{theorem}[\bf Linear Run Time for Ergodic DMDP] \label{thm-2}
Suppose that the Markov decision process specified by $\mathcal{M}=(\cS,\cA,\cP, \br,\gamma)$ is ergodic in the sense that $c_1 \bq\leq \nu^{\bmu} \leq c_2 \bq$ for some distribution vector $\bq$ and any policy $\pi$. 
For any $\epsilon \in(0,1)$, $\delta\in(0,1)$, let $\theta = 1-\gamma + \gamma\frac{c_1}{c_2}$, $ T = \Omega\lf(\lf(\frac{c_2}{c_1}\ri)^4\frac{\S \A \log(\S \A) }{(1-\gamma)^4\epsilon^2}   \ri).$
Then:
\begin{enumerate}[(i)]
\item Algorithm \ref{algo:primaldual} outputs a policy $\hat\bmu$ satisfying $ \bq^{\top }\bv^{\hat \bmu} \geq  \bq^{\top }\bv^*-\epsilon$ with probability at least $2/3$.
\item
Algorithm \ref{algo:meta} outputs a policy $\hat\bmu$ such that $ \bq^{\top }\bv^{\hat \bmu} \geq  \bq^{\top }\bv^*-\epsilon$ in run time
\begin{equation}\label{eq-runtime}
  \tO\lf(\S^2\A +  \lf(\frac{c_2}{c_1} \ri)^4 \frac{\S \A  }{(1-\gamma)^4\epsilon^2} \log\lf(\frac1\delta\ri)  + \frac{1}{\epsilon^2(1-\gamma)^2}\lf(\log\frac1\delta\ri)^2 \ri)
  \end{equation}
 with probability at least $1-\delta$. 
 \end{enumerate}
\end{theorem}

\pf (i)
We first verify that $\mu^*\in\cU$ with the given vector $\bq$ and $\theta= 1-\gamma + \gamma \frac{c_1}{c_2} $.
We note that all eigenvalues of a probability transition matrix $P^{\pi}$ belong to the unit circle and $\gamma\in(0,1)$, therefore the eigenvalues $I-\gamma P^{\pi}$ belong to the positive half plane. As a result, the matrix $I-\gamma P^{\pi}$ is invertible for all $\pi$, including $I-\gamma P^*$.
We have $\suma\mu^*_a
 = \lf(I-\gamma (P^*)^{\top} \ri)^{-1}  (1-\gamma ) \bq $, therefore
\begin{align*}
\suma\mu^*_a
 =  (1-\gamma ) \lf(\sum_{k=0}^{\infty} (\gamma P^*)^k \ri)^{\top} \bq 
&= (1-\gamma) \bq + (1-\gamma ) \lf( \sum_{k=1}^{\infty} (\gamma P^*)^k \ri)^{\top} \bq\\
&\geq (1-\gamma) \bq + \frac1{c_2} (1-\gamma )\sum_{k=1}^{\infty} \lf( (\gamma P^*)^k\ri)^{\top} \nu^* \\
&= (1-\gamma) \bq +  \frac1{c_2}   (1-\gamma )\lf(\sum_{k=1}^{\infty} \gamma^k\ri)\nu^* \\
&\geq (1-\gamma) \bq + \frac{c_1}{c_2} \gamma \bq \\
&\geq  \lf(1-\gamma + \gamma \frac{c_1}{c_2} \ri) \bq. 
\end{align*}
As a result, we have verified $\mu^*\in\cU$ and the assumption and results of Prop.\ \ref{thm-dualitygap} hold.

When Algorithm \ref{algo:primaldual} is applied with $\bq$ and $\theta= 1-\gamma + \gamma \frac{c_1}{c_2} $.
By the nature of the algorithm, we have
$\suma \hat \bpi_a  \geq \theta \bq$. Then we have
and 
$c_1 \bq \geq \nu^{\bmu} \geq c_2 \bq$ for any policy $\bmu$, so the assumptions of Prop.\ \ref{prop-gaptoav} hold. We denote for short that $Gap = \sum_{a\in\cA,i\in\cS}\hat \bpi_{a,i} \lf( v_i^* - \a \sumj p_{ij}(a) v^*_j -\sumj p_{ij}(a) r_{ij}(a) \ri)$. 
By applying Prop.\ \ref{prop-gaptoav}, we obtain that $ \bq^{\top} (\bv^*-\bv^{\hat\pi}) \leq \frac{c_2^2} {c_1^2(1-\gamma )} Gap$. Then we may use the Markov inequality to show that $ \bq^{\top} (\bv^*-\bv^{\hat\pi}) \leq \epsilon$ with probability $2/3$ as long as $ \frac{c_2^2} {c_1^2(1-\gamma )} \E\[ Gap \] \leq 2/3\epsilon$, which requires $ T = \Omega\lf( \lf(\frac{c_2}{c_1}\ri)^4\frac{\S \A \log(\S \A) }{(1-\gamma)^4\epsilon^2}   \ri)$ according to Prop.\ \ref{thm-dualitygap}.

 \vspace{4pt}
 
\noindent (ii) By using a similar analysis as in the proof of Theorem \ref{thm-1}, we finish the proof.
\qed


 \vspace{6pt}
 
Theorem \ref{thm-2} shows that the iteration complexity of Algorithms \ref{algo:primaldual}-\ref{algo:meta} substantially improves when the DMDP is ergodic under every stationary policy. More specifically, the complexity reduces when all policies generate ``similar" stationary distributions. Comparing the preceding complexity result with the input size $\cO(\S^2\A)$ of the DMDP, we conclude that Algorithms \ref{algo:primaldual}-\ref{algo:meta} have nearly-linear run-time complexity in the worst case. For large-scale problems, as long as $\S\gg \frac1{\epsilon}$ and $\A \gg \frac1{\epsilon}$, Algorithm \ref{algo:primaldual} is able to compute an approximate policy more efficiently than any known deterministic method.

The ratio $\frac{c_2}{c_1}$ characterizes a notion of complexity of ergodic DMDP, i.e., the range of possible ergodic distributions under different policies.
Intuitively, dynamic programs are harder to solve when different policies lead to vastly different state trajectories. For example, suppose that there is a critical state that can only show up after a particular sequence of correct actions (this typical happens in aperiodic Markov chains). In this case, any algorithm would need  to search over the space of all policies to identify the critical state. 
For another example, consider that all policies only affect the immediate reward but will lead to the same outgoing transition probabilities. In this case, one would be able learn the optimal action at each state much more efficiently, where $c_1/c_2 =1$.

Our last result shows that it is possible to skip the preprocessing step, as long as the DMDP tuple is specified using special formats. When the input data are given in a way that enables immediate sampling, one can skip Step 5 in Algorithm \ref{algo:primaldual} and remove the first term $\S^2\A$  from the overall run time.


\begin{theorem}[\bf Sublinear Time Complexity for Ergodic DMDP In Special Formats] \label{thm-3}
Suppose that the assumptions of Theorem \ref{thm-2} hold and the collection of transition probabilities $\cP = (P_a)_{a\in\cA}$ are specified in any one of the following formats:
\begin{enumerate}[(a)]
\item Matrices of cumulative probabilities $C_a$ of dimension $\cS\times\cS$, for all $a\in\cA$, where $C_a(i,j) = \sum_{k=1}^j P_a(i,k)$ for all $i \in \cS, a \in \cA$ and $j \in \cS$.
\item Transition probability distributions $P_a(i,\cdot)$ that are encoded in binary trees. There are $\S\A$ trees, one for each state-action pair $(i,a)\in\cS\times\cA$. Each binary tree is of depth $\log|\cS|$ and has $|\cS|$ leaves that store the values of $P_a(i, j)$, $j\in\cS$. Each inner node of the tree stores the sum of its two children. 
\end{enumerate}
Then for any $\epsilon \in(0,1)$, $\delta\in(0,1)$,
Algorithm \ref{algo:meta} outputs an approximately optimal policy $\hat\bmu$ such that 
$ \bq^{\top }\bv^{\hat \bmu} \geq  \bq^{\top }\bv^* -\epsilon$ 
in run time
$$  \tO\lf(  \lf(\frac{c_2}{c_1} \ri)^4 \frac{\S \A  }{(1-\gamma)^4\epsilon^2} \log\lf(\frac1\delta\ri)  + \frac{1}{\epsilon^2(1-\gamma)^2}\lf(\log\frac1\delta\ri)^2 \ri)$$
 with probability at least $1-\delta$. 
\end{theorem}
\pf  In both cases (a) and (b), it is possible to  draw a sample coordinate $j\mid (i,a)$ with probability $p_{ij}(a)$ using $\cO(\log |\cS|)$ run time, because the corresponding data structure storing the vector $P_a(i,\cdot)$ can be readily used as a sampler  \cite{bringmann2012efficient, wong1980efficient}. Therefore one can skip the preprocessing step (Step 5) in Algorithm \ref{algo:primaldual} and remove the $\tO(\S^2\A)$ term from the run-time result of Theorem \ref{thm-2}.
\qed

\vspace{6pt}

When the input is given in suitable data structures, the run-time complexity of Algorithms \ref{algo:primaldual}-\ref{algo:meta} reduces from nearly-linear to sublinear with respect to the input size.
The reduced run time is almost linear in $|\cS\times \cA| $, i.e., the number of state-action pairs. It means that for fixed values of $\epsilon,\gamma$, each state-action pair is queried for a constant number of times on average regardless of the dimension of the DMDP. This result is sublinear with respect to the input size $\cO(\S^2\A).$ It suggests that one can find an approximately optimal policy without even reading a significant portion of the input entries. We recall that the input data mainly consist of transition probabilities $p_{ij}(a)$. An explanation for the sublinear complexity is that certain small probabilities in the input data can be safely ignored, without deteriorating the quality of the approximate policy significantly.

Theorems \ref{thm-1}, \ref{thm-2}, \ref{thm-3} established new complexity upper bounds for computing an approximate-optimal policy of the DMDP. Let us compare the upper bounds given by Theorems \ref{thm-1}, \ref{thm-2}, \ref{thm-3} with recent lower bound results for DMDP \cite{chen2017lower}. For DMDP that is specified in the standard way (arrays of transition probabilities), it shows that any randomized algorithm needs at least $\Omega(\S^2 \A)$ run time to get any $\epsilon$-approximate policy with high probability; see Theorem 1 of \cite{chen2017lower}. It also shows that the lower bound reduces to  $\Omega(\frac{\S\A}{\epsilon})$ when the input data are in the format of binary trees or cumulative sums (for which the preprocessing step can be skipped); see Theorems 2-3 of \cite{chen2017lower}. Our upper bound results nearly match the lower bound results in both cases. Both the upper and lower bounds suggest that the computational complexity of DMDP indeed depends on the input data structure. The overall complexity for approximately solving the MDP is dominated by the preprocessing time. Once the data is preprocessed, the remaining computation problem becomes significantly easier.

\section{Remarks} 

We have developed a novel randomized method that exploits the linear duality between the value function and the policy function for solving DMDPs. 
It is related to several fundamental methods in linear programming, stochastic optimization and online learning. It is easy to implement, uses sublinear space complexity and nearly-linear (sometimes sublinear) run-time complexity.

Algorithm \ref{algo:primaldual} can be viewed as a randomized version of the simplex method, which is also equivalent to a version of the policy iteration method for solving MDP. It maintains a primal variable (value) and a dual variable (policy). Each update $\bmu_{i,a} \leftarrow \bmu_{i,a} \cdot \exp\lf\{\g \ri\}$ mimics a pivoting step towards a neighboring basis. Instead of moving from one basis to another one, each update in the dual variable can be viewed as a ``soft" pivot. Our theoretical results show that the algorithm finds an approximate policy in $\tO(SA)$ iterations. This is somewhat similar to the simplex method, which also terminates in $\tO(SA)$ (which is the number of constraints) iterations on average. What makes our randomized algorithm different is its $\tO(1)$ run time per iteration. It avoids explicitly solving any linear system, which is unavoidable in the simplex method.

Algorithm \ref{algo:primaldual} can also be viewed as a stochastic approximation method for solving a saddle point problem. It utilizes the structures of specially crafted primal and dual constraint sets to make each update as simple and efficient as possible. The update of the dual variable (policy) uses a special Bregman divergence function which is related to the relative entropy between randomized policies. 

It is worth noting that Algorithm \ref{algo:primaldual} is related to the exponentiated gradient method for the multi-arm bandit problem in the online learning setting. When there is a single state, Algorithm \ref{algo:primaldual} reduces to the basic exponentiated gradient method. This observation provides a hint that we might be able to adapt Algorithm \ref{algo:primaldual} to apply to online reinforcement learning. This is a direction for future research.

The new method and complexity results of this paper suggest a promising direction that awaits further research. The current results leave open many questions. We conjecture that the complexity result should have a better dependence on $\frac1{1-\gamma}$, especially when all $P^{\pi}$'s have relatively large spectral gaps. A related notion of complexity metric for MDP is the {diameter}, i.e., the maximal expected time to move from any state to any other state. We conjecture that the diameter should play a key role in an improved complexity analysis and replace at least one factor of $\frac1{1-\gamma}$.
We also conjecture that the sublinear-time complexity result of Proposition 3 hold for more general DMDPs without prior knowledge about $\frac{c_2}{c_1}$. 
Another direction for future research is to consider the run-time complexity for finite-horizon MDP and average-reward MDP. It remains unclear what roles the mixing rate and the horizon play in the run-time complexity. In the mean time, an equally important (if not more) question is to establish the computation complexity lower bound for approximating optimal policies of MDP.

\section{Proof of Proposition \ref{thm-dualitygap}}

We provide the complete proof of Proposition \ref{thm-dualitygap} in this last section. Readers who are not interested in the technical details are free to skip this part. 

\subsection{Technical Lemmas}\label{sec-lemmas}
 
In this section, we analyze the convergence of Algorithm \ref{algo:primaldual}.
We denote the $t$-th iterates generated by Algorithm \ref{algo:primaldual} by $\pi^t$, $\xi^t$, and $\bv^t$.
We define the auxilary variables $\lambda = \lf(\lambda^t_{i,a}\ri)_{i\in\cS,a\in\cA}$ as 
$$\lambda^t_{i,a} = \xi^t_i \pi^t_{i,a}.$$
According to Algorithm \ref{algo:primaldual}, we can verify that $\xi^t\in\Re^{\S}$ and $\pi^t_i\in\Re^{\A}$ are vectors of probability distributions. It follows that $\lambda^t$ is always a $\S\A$-dimensional vector of probability distribution. 
In addition,  
the updates on $\xi^t$ and $\pi^t$ can be equivalently written as updates on $\lambda^t$ and $\pi^t$, given by
\begin{equation}\label{eq-update}
\lambda^{t+1}_{{i,a}} = \frac{\lambda^{t}_{{i,a}} \cdot \exp(\g^{t+1}_{{i,a}})}{\sum_{i',a'}
\lambda^{t}_{i',a'} \cdot \exp(\g^{t+1}_{i',a'})}, \qquad
\pi^{t+1}_{{i,a}} = \frac{\pi^{t}_{{i,a}} \cdot \exp(\g^{t+1}_{{i,a}})}{\sum_{a'}
\pi^{t}_{i,a'} \cdot \exp(\g^{t+1}_{i,a'})}
\qquad \forall\ i\in\cS,a\in\cA,
\end{equation}
where
\begin{equation}\label{eq-gg}
\g^{t+1}_{{i,a}} =
\lf\{
\begin{tabular}{l l}
$ {\beta}  \cdot \frac{(\gamma v^t_j - v^t_{i} +r_{ij^t}(a) -M )}{ ((1-\theta) \xi_i^t +\theta q_i)\pi^t_{{i,a}}}$
 & if $i=i_{t+1}, a = a_{t+1}$
\\
0  & otherwise
\end{tabular}
\ri.
\end{equation}
In what follows, we denote by $\cF_t$ the collection of random variables that are revealed up to the end of the $t$-th iteration.
We denote by $\mu^t_{i,a}$ the dual iterates given by
$$\mu^t_{i,a} = ((1-\theta) \xi_i^t +\theta q_i)\pi_{i,a}^t = (1-\theta)\lambda_{i,a}^t + \theta q_i \pi_{i,a}^t.$$
According to Algorithm \ref{algo:primaldual}, we can verify that $\mu^t\in\cU$ and $\bv^t\in\cV$ for all $t$ with probability 1.

\begin{lemma}
If $\mu^*\in\cU$, there exists probability distribution vectors $\lambda^* \in\Re^{\S\A}$ and $\pi_i^* \in\Re^{\A}$, $i\in\cS$, such that
$$\bpi^*_{i,a} =  (1-\theta)  \lambda^*_{i,a} + \theta  q_i \bmu^*_{i,a},\qquad \forall~i\in\cS,a\in\cA.$$
\end{lemma}
\pf The proof is straightforward by the definition of $\cU$.
\qed

\begin{lemma} \label{lemma-KL}The iterates generated by Algorithm \ref{algo:primaldual} satisfy
\begin{equation}\label{eq-KL0}
\begin{split}
& \E\[\KL(\lambda^* || \lambda^{t+1}) \mid \cF_t\] - \KL(\lambda^* || \lambda^t)  
\leq \sumia (\lambda^t_{{i,a}}  - \lambda^*_{{i,a}}) \E\[  \g^{t+1}_{{i,a}} \mid \cF_t\]+  \frac12\sumia \lambda^t_{{i,a}}\E\[  \lf(\g^{t+1}_{{i,a}}\ri)^2\mid \cF_t\] ,
\end{split}
\end{equation}
for all $t$, with probability 1.
\end{lemma}

\pf
\def\bg{\mathbf{g}}
By using the relation \eqref{eq-update}, we have
\begin{equation}\label{eq-01}\begin{split}
\KL(\lambda^* || \lambda^{t+1}) - \KL(\lambda^* || \lambda^t) 
&= \sumia \lambda^*_{{i,a}} \log \frac{\lambda^*_{{i,a}}}{\lambda^{t+1}_{{i,a}}} - \sumia \lambda^*_{{i,a}} \log \frac{\lambda^*_{{i,a}}}{\lambda^t_{{i,a}}} \\
&= \sumia \lambda^*_{{i,a}} \log \frac{\lambda^t_{{i,a}}} {\lambda^{t+1}_{{i,a}}} \\
&= \sumia \lambda^*_{{i,a}} \log \frac{Z } {\exp( \g^{t+1}_{{i,a}})} \\
&= \sumia \lambda^*_{{i,a}} \log \lf(Z\ri)  - \sumia \lambda^*_{{i,a}} \g^{t+1}_{{i,a}}\\
&=  \log Z  -\sumia \lambda^*_{{i,a}} \g^{t+1}_{{i,a}},
\end{split}\end{equation}
where $Z = \sumia \lambda^t_{{i,a}} \exp(\g^{t+1}_{{i,a}} )$.
According to \eqref{eq-gg}, we have $\gamma v^t_j - v^t_{i} +r_{ij^t}(a) -M \leq \frac{\gamma}{1-\gamma}-0+1 -\frac1{1-\gamma}\leq 0$ because $v_i\in[0,\frac1{1-\gamma}]$, $r_{ij^t}(a)\in[\frac12,1]$ and $M=\frac1{1-\gamma}$. It follows that $ \g^{t+1}_{i,a} \leq 0$ for all $i\in\cS,a\in\cA$ with probability 1. Then we derive
\begin{equation}\label{eq-02}\begin{split}
\log Z =  \log \lf( \sumia \lambda^t_{{i,a}} \exp(\g^{t+1}_{{i,a}} )\ri) 
&\leq \log \sumia \lambda^t_{{i,a}} \lf(1+ \g^{t+1}_{{i,a}} + \frac12 \lf(\g^{t+1}_{{i,a}}\ri)^2 \ri)\\
&= \log \lf( 1 + \sumia \lambda^t_{{i,a}}  \g^{t+1}_{{i,a}} +\frac12 \sumia \lambda^t_{{i,a}}  \lf(\g^{t+1}_{{i,a}}\ri)^2 \ri)\\
&\leq \sumia \lambda^t_{{i,a}}  \g^{t+1}_{{i,a}} + \frac12 \sumia \lambda^t_{{i,a}}  \lf(\g^{t+1}_{{i,a}}\ri)^2 ,
\end{split}\end{equation}
where the first inequality uses the fact $e^x\leq 1+x+\frac12 x^2$ if $x\leq 0$ and the second inequality uses the fact $\log(1+x)\leq x$ for all $x$.
We combine \eqref{eq-01} and \eqref{eq-02} and take conditional expectation $\E\[\cdot\mid \cF_t\]$ on both sides, then we obtain \eqref{eq-KL0}.
\qed

\begin{lemma} \label{lemma-KL-pi}
For any $i\in\cS$,
the iterates generated by Algorithm \ref{algo:primaldual} satisfy
\begin{equation}\label{eq-KL-pi}
\begin{split}
\E\[\KL(\bmu^*_i || \bmu^{t+1}_i) ~\Big\vert~ \cF_t\] - \KL(\bmu^*_i || \bmu_i^t) 
&\leq \suma (\bmu^t_{i,a} - \bmu^*_{i,a}) \E\[ \g^{t+1}_{i,a} \mid \cF_t\]+ \frac12 \suma \bmu^t_{i,a}  \E\[\lf(\g^{t+1}_{i,a}\ri)^2\mid \cF_t\] ,
\end{split}
\end{equation}
for all $t\geq 1$, with probability 1.
\end{lemma}
\pf
The proof is similar to Lemma \ref{lemma-KL}. We omit it for simplicity.
\qed

\vspace{5pt}

\begin{lemma}\label{lemma-variance}
The iterates generated by Algorithm \ref{algo:primaldual} satisfy 
$$\sumia \bpi^t_{{i,a}}\E\[  \lf(\g^{t+1}_{{i,a}}\ri)^2\mid \cF_t \] \leq \frac{4 \S\A\beta^2 }{(1-\gamma)^2 },$$
for all $t\geq 1$ with probability 1.
\end{lemma}
\pf We note that $\mathbf{P}(i_{t+1} =i,a_{t+1} = a\mid \cF_t) = \mu^t_{i,a}.$ Then we have
\begin{align*}
\sumia  \bpi^t_{{i,a}} \E\[  \lf( \g^{t+1}_{{i,a}}\ri)^2\mid \cF_t \]
&= \sumi\suma  \bpi^t_{{i,a}} \cdot  \bpi^t_{{i,a}} \cdot \sumj p_{ij}(a)\lf( {\beta}  \cdot \frac{(\gamma v^t_j - v^t_{i} +r_{ij}(a)   -M )}{ \bpi^t_{{i,a}}}\ri)^2
\\
&= \sumi\suma \sumj p_{ij}(a) \lf( {\beta}   \cdot  (\gamma v^t_j - v^t_{i} +r_{ij}(a) -M ) \ri)^2\\
&\leq \sumi\suma  \lf(\beta \cdot  \frac{2 }{1-\gamma} \ri)^2\\
&=\frac{ 4\S\A  \beta^2}{(1-\gamma)^2 },
\end{align*}
where the inequality uses the fact that $\bv^t$ belongs to the $\|\cdot\|_{\infty}$ ball of radius $\frac{ 1 }{1-\gamma}$ and $M=\frac1{1-\gamma}$.
\qed




\begin{proposition} 
We let $\Phi(\mu^t )$ be the divergence function given by
$$\Phi(\mu^t ) = (1-\theta) \KL(\lambda^* || \lambda^t) + \theta \sumi   q_i \KL(\bmu^*_i || \bmu_i^t) .$$
The iterates generated by Algorithm \ref{algo:primaldual} satisfy
\begin{equation}\label{eq-Phi}
\E\[ \Phi(\mu^{t+1} ) \mid \cF_t\]  \leq \Phi(\mu^t ) \\+\beta \suma (\bpi_a^t - \bpi_a^*)^{\top} \lf( ( \gamma P_a - I ) \bv^t + \br_a\ri) +  \frac{2\S\A \beta^2}{(1-\gamma)^2 },
\end{equation}
for all $t$, with probability 1.
\end{proposition}

\pf We take the weighted sum between \eqref{eq-KL0} and \eqref{eq-KL-pi}, so we have
\begin{equation*}
\begin{split}
\E\[ \Phi(\mu^{t+1} ) \mid \cF_t\] 
&\leq \Phi(\mu^t ) + (1-\theta) \sumia  (\lambda^t_{{i,a}}  - \lambda^*_{{i,a}}) \E\[  \g^{t+1}_{{i,a}} \mid \cF_t\]+ \frac{1-\theta}2 \sumia\lambda^t_{{i,a}}  \E\[  \lf(\g^{t+1}_{{i,a}}\ri)^2\mid \cF_t\]
\\
&\qquad + \theta \sumi  q_i \suma  (\bmu^t_{{i,a}}  - \bmu^*_{{i,a}}) \E\[  \g^{t+1}_{{i,a}} \mid \cF_t\]
+ \frac{\theta}2\sumi q_i\suma \pi^t_{{i,a}} \E\[  \lf(\g^{t+1}_{{i,a}}\ri)^2\mid \cF_t\]
\\
&= \Phi(\mu^t ) + \sumia (\bpi^t_{{i,a}}  - \bpi^*_{{i,a}}) \E\[  \g^{t+1}_{{i,a}} \mid \cF_t\]+  \frac12\sumia \bpi^t_{{i,a}}\E\[  \lf( \g^{t+1}_{{i,a}}\ri)^2\mid \cF_t\]
,
\end{split}\end{equation*}
where the equality uses the relations $\bpi^t_{i,a} = (1-\theta)  \lambda^t_{i,a} + \theta  q_i \bmu^t_{i,a}$
and $\bpi^*_{i,a} =  (1-\theta)  \lambda^*_{i,a} + \theta  q_i \bmu^*_{i,a}$ (by Lemma 1 since $\mu^*\in  \cU$).
For arbitrary $i\in\cS$ and $a\in\cA$, we have
$$\frac{1}\beta\cdot \E\[ \g^{t+1}_{{i,a}}  \mid \cF_t\] =  \gamma\sumj p_{ij}(a)v^t_j  -v^t_{i} + \sumj p_{ij}(a) r_{ij}(a) -M  = (\gamma P_a \bv^t - \bv^t + \br_a)_i - M.$$
It follows that
\begin{align*}
\frac{1}\beta\cdot \sumia (\bpi^t_{{i,a}}  - \bpi^*_{{i,a}}) \E\[ \g^{t+1}_{{i,a}}  \mid \cF_t\] 
&= \suma\sumi(\bpi^t_{i,a} - \bpi^*_{i,a})\[ 
(\gamma P_a \bv^t - \bv^t + \br_a)_i - M
 \] \\
 &= \suma (\bpi_a^t - \bpi_a^*)^{\top} \lf( ( \gamma P_a - I ) \bv^t + \br_a\ri),
\end{align*}
where the second equality comes from the fact $\sumia \bpi^t_{{i,a}} = \sumia  \bpi^*_{{i,a}} =1$ (because $\bpi^t\in\cU$, $\bpi^*\in\cU$).
According to Lemma \ref{lemma-variance}, we also have
$\E\[\sumia \bpi^t_{{i,a}}  \lf(\g^{t+1}_{{i,a}}\ri)^2\mid \cF_t\] \leq  \frac{4\S\A \beta^2}{(1-\gamma)^2 } .$
We apply the preceding two relations and complete the proof.

\qed



\begin{proposition}\label{lemma-v}
The iterates generated by Algorithm \ref{algo:primaldual} satisfy for all $t$ with probability 1 that
\begin{equation}\label{eq-v}
\E\[\|\bv^{t+1}-\bv^*\|^2 \mid \cF_t\]
\leq \|\bv^t-\bv^*\|^2 + 2 {\alpha} (\bv^t-\bv^*)^{\top} 
\lf(\suma  (I-\gamma P_a)^{\top}\bpi_a^t - (1-\gamma) \bq \ri)  + {4\alpha^2} .
\end{equation}
\end{proposition}

\pf 
We let $\bv^{t+1/2} $ be the vector such that
\begin{align*}
v^{t+1/2}_{i_{t+1}} & =
v^t_{i_{t+1}} - {\alpha}\lf( \frac{(1-\gamma)q_{i_{t+1}} }{(1-\theta) \xi^t_{i_{t+1}} +\theta q_{i_{t+1}}} -1 \ri)
,\\
v^{t+1/2}_{j_{t+1}} & =
v^t_{j_{t+1}} -  {\alpha} \gamma 
,\\
v^{t+1/2}_{ i} &= v^t_{i}, \qquad \hbox{if   } i\notin\{ i_{t+1}, j_{t+1}\}.
\end{align*}
Then we can verify that
$\bv^{t+1} = \Pi_{\cV}\bv^{t+1/2} ,$ where $\Pi$ denotes the Euclidean projection.
We note that $\P(i_{t+1} = i \mid \cF_t) =(1-\theta) \xi^t_{i} +\theta q_{i} = \suma \mu^t_{i,a}$,
$\P(j_{t+1} =j \mid \cF_t) = \sumia \mu^t_{i,a} p_{ij}(a)$.
Then we can verify that
$$\E\[\bv^{t+1/2}-\bv^t \mid \cF_t\] =   -{ \alpha} \lf((1-\gamma) \bq - \suma  (I-\gamma P_a)^{\top}\bpi_a^t \ri).
$$
Also since $\theta\geq 1-\gamma$ we have $\frac{(1-\gamma)q_{i_{t+1}}}{ (1-\theta)\xi^t_{i_{t+1}} + \theta q_{i_{t+1}}}\in[0,1]$. Then we have  $|v^{t+1/2}_{i_{t+1}} - v^{t}_{i_{t+1}} | < {\alpha}  $ and $|v^{t+1/2}_{j_{t+1}} - v^{t}_{j_{t+1}}|< {\alpha} $. Then we can verify that
$\|\bv^{t+1/2}-\bv^t\|^2 \leq {4\alpha^2}$
for all $t$ with probability 1. 
Finally, by using the nonexpansive property of $\Pi_{\cV}$ and $\bv^*\in \cV$, we further obtain
\begin{align*}
\E\[\|\bv^{t+1}-\bv^*\|^2 \mid\cF_t\]
&= \E\[\|\Pi_{\cV}\bv^{t+1/2}-\bv^*\|^2 \mid\cF_t\]\\
&\leq \E\[\|\bv^{t+1/2}-\bv^*\|^2 \mid\cF_t\]\\
&= \|\bv^t-\bv^*\|^2 +2(\bv^t-\bv^*)^{\top}\E\[(\bv^{t+1/2}-\bv^t)\mid\cF_t\]  + \E\[\|\bv^{t+1/2}-\bv^t\|^2\mid\cF_t\],
\\&\leq\  \|\bv^t-\bv^*\|^2 
+ {2\alpha} (\bv^t-\bv^*)^{\top} \lf(\suma  (I-\gamma P_a)^{\top}\bpi_a^t - (1-\gamma) \bq \ri)  
+{4\alpha^2},
\end{align*}
for all $t$ with probability 1.
\qed


%
%

\begin{proposition}
We define for short that 
$$\cE^t = \Phi(\mu^t) + \frac{(1-\gamma)^2}{\S } \|\bv^t-\bv^*\|^2$$
and
$$\cG^t  =  \suma (\bpi_a^t)^{\top} \lf( (I- \gamma P_a  ) \bv^* - \br_a\ri)  =  \sumia \bpi^t_{{i,a}}( \bv^* - \a P_{a} \bv^* - \br_{a} )_i. $$
Let $\alpha = \frac{\S}{2(1-\gamma)^2} \beta$.
The iterates generated by Algorithm \ref{algo:primaldual} satisfy for all $t$ with probability 1 that
\begin{equation}\label{thm1:equ5}
\E\[\cE^{t+1}\mid \cF_t\]
\leq \cE^{t} -  \beta  \cG^t 
+ {\beta^2}    \frac{2\S(\A +1) }{(1-\gamma)^2}.
\end{equation}
\end{proposition}

\pf Let $\alpha = \frac{\S}{2(1-\gamma)^2} \beta$. We multiply \eqref{eq-v} with $\frac{(1-\gamma)^2}{\S }$ and takes its sum with \eqref{eq-Phi}, obtaining
\begin{equation*}
\begin{split}
\E\[\cE^{t+1}\mid \cF_t\]
&\leq \cE^{t} + {\beta^2}  \frac{2\S\A  +\S}{(1-\gamma)^2}  \\
&+ {\beta}  \lf( \suma (\bpi_a^t - \bpi_a^*)^{\top} \lf( ( \gamma P_a - I ) \bv^t + \br_a\ri)  + (\bv^t-\bv^*)^{\top} \lf(\suma  (I-\gamma P_a)^{\top}\bpi_a^t - (1-\gamma) \bq \ri) \ri) .
\end{split}
\end{equation*}
We have
\begin{align*}
& \suma (\bpi_a^t - \bpi_a^*)^{\top} \lf( ( \gamma P_a - I ) \bv^t + \br_a\ri)  + (\bv^t-\bv^*)^{\top} \lf(\suma  (I-\gamma P_a)^{\top}\bpi_a^t - (1-\gamma) \bq \ri)
\\
&= \suma (\bpi_a^t - \bpi_a^*)^{\top} \lf( ( \gamma P_a - I ) \bv^t + \br_a\ri)  + (\bv^t-\bv^*)^{\top} \suma  (I-\gamma P_a)^{\top} ( \bpi_a^t -  \bpi_a^*)\\
&= \suma (\bpi_a^t - \bpi_a^*)^{\top} \lf( ( \gamma P_a - I ) \bv^* + \br_a\ri) \qquad (\hbox{by the dual feasibility of $\mu^*$}) \\
&= \suma (\bpi_a^t)^{\top} \lf( ( \gamma P_a - I ) \bv^* + \br_a\ri) \qquad (\hbox{by the linear complementarity condition for $\bv^*$, $\mu^*$})  
 \end{align*}
 where the second equality uses the dual feasibility of $\mu^*$ that $\suma (I-\gamma P_a)^{\top}\bpi^*_a = (1-\gamma)\bq$ and the fourth equality uses the complementary condition
 $\bpi_{a,i}^* \lf( ( \gamma P_a - I ) \bv^* + \br_a\ri)_i =0$ for all $i\in\cS,a\in\cA$.
 Combining the preceding relations, we obtain  \eqref{thm1:equ5}. 
\qed

\subsection{Proof of Proposition \ref{thm-dualitygap}} 

\pf
We claim that $\cE^1\leq \log(\S\A) + 1$. To see this, we note that $\lambda^1$ and $\pi^1_i$'s are uniform distributions (according to Step 3 of Algorithm \ref{algo:primaldual}). Therefore we have $\KL(\lambda^*||\lambda^1)\leq \log (SA)$ and $\KL(\bmu_i^*||\bmu_i^1)\leq \log (S)$ for $i$, and $\|\bv^t- \bv^*\|^2\leq \frac{S}{(1-\gamma)^2}$ for all $t$. Then we have $\cE^1 \leq (1-\theta) \KL(\lambda^* || \lambda^1) + \theta \sumi   q_i \KL(\bmu^*_i || \bmu^1_i)  + \frac{(1-\gamma)^2}{\S } \|\bv^1-\bv^*\|^2 \leq \log(\S\A) + 1$.

We rearrange the terms of \eqref{thm1:equ5} and obtain
$$
 \gk 
      \leq \frac{1}{\beta }(\ek-  \E\[ \cE^{t+1} \mid \cF_t\])+\frac{2\beta \S(\A+1) }{(1-\gamma)^2}.
$$
Summing over $t=1,\ldots,T$ and taking average, we have 
   \begin{equation*}
       \begin{aligned}
             \E\[  \sum_{t=1}^{T}  \gk \] 
           &\leq \frac{1}{ \beta  }\sum_{t=1}^{T}( \E\[\ek\]- \E\[ \cE^{t+1} \])+\frac{2\beta \S(\A+1) T}{(1-\gamma)^2}\\
           &= \frac{\E\[\cE^1\]-\E\[\cE^T\]}{\beta } + \frac{2\beta \S(\A+1) T}{(1-\gamma)^2}\\
           &\leq \frac{1}{\beta  } (\log(\S\A)+1)+ \frac{2\beta \S(\A+1)T}{(1-\gamma)^2}.
     \end{aligned}
   \end{equation*}
where the inequality is based on the fact $\cE^1 \leq \log(\S\A) + 1$ and $\cE^T\geq 0$. 
Therefore by taking $\beta  =(1-\gamma) \sqrt{\frac{\log \S\A +1}{2\S\A T} }$,
we obtain
$
  \E\[  \frac{1}{T}\sum_{t=1}^{\top}\gk \]   \leq \frac{\sqrt{ 2\S(\A+1)(\log \S\A+1)} }{(1-\gamma)\sqrt{T}} .
$
\qed

\def\bibfont{\footnotesize}
 \bibliography{nips_1705,mdp,mdp2,mdp3}
{\bibliographystyle{plain}}

\end{document}